\documentclass[12pt]{article}
\usepackage{ct}

\usepackage{mathtools,tikz,longtable,multirow}

\theoremstyle{plain}
\newtheorem{const}[theorem]{Construction}

\specs{n}{vol (issue)}{year}{}

\dateline{Jan 1, 2020}{Jan 2, 2020}{TBD}



\title{New Constructions of Distance-Biregular Graphs}

\author[1]{Blas Fern\'andez\thanks{Supported by the Slovenian Research and Innovation Agency (ARIS) through the research program P1-0285 and the research projects J1-50000 and Z1-70008.}}
\author[2]{Ferdinand Ihringer}
\author[3]{Sabrina Lato\thanks{Supported by Kempe foundation JCSMK22-0160.}}
\author[4,5]{Akihiro Munemasa\thanks{Supported by JSPS KAKENHI Grant Number 25K07095.}}

\affil[1]{%
Inštitut za Matematiko, Fiziko in Mehaniko, Ljubljana, Slovenia and Faculty of Mathematics, Natural Sciences, and Information Technologies, University of Primorska, Koper, Slovenia.

\email{blas.fernandez@famnit.upr.si}
}

\affil[2]{%
Dept.~of Mathematics,
Southern University of Science and Technology, Shenzhen, Guangdong, China.

\email{ihringer@sustech.edu.cn}
}

\affil[3]{%
  Faculty of Mathematics, Natural Sciences, and Information Technologies,
  University of Primorska, Koper, Slovenia.

\email{sabrina.lato@iam.upr.si}
}

\affil[4]{%
Graduate School of Information Sciences,
Tohoku University, Sendai, Japan.

\email{munemasa@tohoku.ac.jp}
}

\affil[5]{%
School of Science, China University of Geosciences, Beijing, China.
}

\newcommand\pmat[1]{\begin{pmatrix} #1 \end{pmatrix}}
\newcommand\vmat[1]{\begin{vmatrix} #1 \end{vmatrix}}

\newcommand{\abs}[1]{\left | #1 \right |}
\newcommand{\lp}{\left (}
\newcommand{\rp}{\right )} 
   \newcommand{\lsb}{\left \{ }
\newcommand{\rsb}{\right \} }

\newcommand{\ver}[1]{\mathbf{e}_{ #1}}
\newcommand\ip[2]{\langle#1,#2\rangle}

\newcommand{\vtxa}{u}
\newcommand{\vtxb}{v}
\newcommand{\vtxc}{w}

\newcommand{\bipartB}{Y}
\newcommand{\bipartC}{Z}
\newcommand{\bipartS}{X}

\newcommand{\qbino}[2]{\left [ #1 \right ]_{#2}}

\newcommand{\<}{\langle}
\renewcommand{\>}{\rangle} 
\newcommand{\qbin}[2]{\genfrac{[}{]}{0pt}{}{#1}{#2}}
\newcommand{\qbinom}[2]{\genfrac{[}{]}{0pt}{}{#1}{#2}}

\newcommand{\cR}{\mathcal{R}}
\newcommand{\cS}{\mathcal{S}}

\newcommand\mysqueeze[2]{%
\newdimen\origspacing%
\newdimen\newspacing%
\origspacing=\fontdimen2\font%
\setlength{\newspacing}{1\origspacing}%
\addtolength{\newspacing}{-#1}%
\fontdimen2\font=\newspacing%
{#2}%
\fontdimen2\font=\origspacing}

\DeclarePairedDelimiter{\norm}{\lVert}{\rVert}

\definecolor{cExist}{RGB}{176,255,176}
\definecolor{cNot}{RGB}{255,176,176}
\definecolor{cNew}{RGB}{176,176,255}
\definecolor{cUnknown}{RGB}{255,255,176}

\begin{document}

\maketitle


\begin{abstract}
 We construct a new family of distance-biregular graphs related to hyperovals
 and a new sporadic example of a distance-biregular graph related to Mathon's perp system.
 The infinite family can be explained using 2-$\bipartB$-homogeneity, while the
 sporadic example belongs to a generalization of a construction by Delorme.
 Additionally, we establish a new non-existence condition for distance-biregular graphs
 which, for instance, rules out the existence of a distance-biregular graph on $225+60$ vertices.
\end{abstract}



\section{Introduction}

Distance-biregular graphs are a class of bipartite graphs with strong algebraic and combinatorial properties. Although the theory of distance-biregular graphs has developed over the years, the known examples have remained essentially stable since they were defined. In this paper, we change that by describing new constructions for an infinite family of distance-biregular graphs, as well as a new sporadic example.

A bipartite graph is \textit{distance-biregular} if for all vertices \( \vtxa \) and \( \vtxb \), the number of vertices adjacent to \( \vtxb \) and closer to \( \vtxa \) than \( \vtxb \) depends only on the distance between \( \vtxa \) and \( \vtxb \) and the cell of the bipartition that \( \vtxa \) lies in. Delorme~\cite{delormeFrench} defined such graphs as having the property ``regularit\'e metrique fort.'' Godsil and Shawe-Taylor~\cite{distanceRegularised} independently studied the same class of graphs under the name ``distance-biregular graphs.'' Delorme later wrote a paper in English~\cite{delorme} which is essentially a translation of his earlier French paper, and we will cite the English version for the relevant results that appear in both.

In addition to the foundational papers of Delorme~\cite{delorme} and Godsil and Shawe-Taylor \cite{distanceRegularised}, early results on distance-biregular graphs can be found in the paper of Mohar and Shawe-Taylor~\cite{biregularCage} or the theses of Van den Akker~\cite{van1990distance} or Shawe-Taylor~\cite{shaweTaylor}. The theory of distance-biregular graphs has also been developed more recently, and the key results that will be relevant for this paper can be found in Fern\'andez and Penji\'c~\cite{fernandez2023almost} and Lato~\cite{phD, soloCharacterization}.

Regular distance-biregular graphs are equivalent to bipartite distance-regular graphs. Dis-tance-regular graphs are a well-studied class of graphs and more information can be found in the monograph of Brouwer, Cohen, and Neumaier~\cite{bcn} or the more recent survey of Van Dam, Koolen, and Tanaka~\cite{vanDistance}. This paper is mostly concerned with graphs that are not distance-regular, so we will generally assume we are dealing with semiregular bipartite graphs with valencies $k \neq \ell.$

Although distance-biregular graphs have been considerably less studied compared to their distance-regular counterparts, they do include several notable families that have been studied for their connections to design theory, finite geometry, and algebra.

A motivating family of distance-biregular graphs that are not distance-regular come from generalized polygons. Tits~\cite{tits1959trialite} introduced generalized polygons, bipartite graphs with girth twice the diameter. A \textit{thick} generalized polygon is a generalized polygon with minimum degree at least three, and any thick generalized polygons are distance-biregular~\cite{delorme, distanceRegularised}. Thick generalized polygons can only exist with diameter \( d =2, 3, 4, 6\),  or \( 8 \)~\cite{generalizedPolyThick}. Infinite families exist for each of these possible diameters, and the only known generalized polygons of diameter eight are biregular but not regular~\cite{van2012generalized}.

Distance-regular graphs of small diameter have been the subject of particular attention. Cvetkovi\'{c}, Doob, and Sachs~\cite{cvetkovic} characterized bipartite distance-regular graphs of diameter three as the incidence graphs of symmetric designs. This similarly characterizes distance-biregular graphs of diameter three, and can be extended to distance-biregular graphs with diameter four and vertices of eccentricity three. Delorme~\cite{delorme} and Shawe-Taylor~\cite{shaweTaylor} proved that such distance-biregular graphs are equivalent to a particular class of quasi-symmetric designs. These quasi-symmetric designs include Steiner systems and affine resolvable designs, and they were studied further from the perspective of distance-biregular graphs in Chapter 5 of Shawe-Taylor~\cite{shaweTaylor}.

Another notable class of distance-regular graphs is the class of distance-regular graphs of diameter two, which are called \textit{strongly regular graphs}. Strongly regular graphs are older than distance-regular graphs, and can be defined more directly. Following the notation of Brouwer and Van Maldeghem~\cite{brouwer2022strongly}, a \textit{strongly regular graph} with parameters $(v, k, \lambda, \mu)$ is a graph on \( v \) vertices with valency \( 1 \leq k \leq v-2 \) such that any two adjacent vertices have $\lambda$ common neighbours and any two nonadjacent vertices have $\mu$ common neighbours. More information can be found in the monograph of Brouwer and Van Maldeghem~\cite{brouwer2022strongly}.

Bose~\cite{bose1963strongly} introduced the notion of partial geometries to study strongly regular graphs. Although partial geometries were defined geometrically, they can equivalently be thought of as distance-biregular graphs with diameter four and girth six. Infinite families of partial geometries exist, and more information can be found in the surveys of Brouwer and Van Lint~\cite{brouwer1984strongly}, De Clerck and Van Maldeghem~\cite{de1995some}, or Thas~\cite{thas2007partial}.

Bipartite distance-regular graphs, generalized polygons, quasi-symmetric designs, and partial geometries are major classes of distance-biregular graphs that are studied in other contexts. Thus, the smallest ``genuinely'' distance-biregular graphs which do not belong to these classes are the non-regular distance-biregular graphs with girth four where every vertex has eccentricity four. However, such graphs still belong to a larger category of graphs studied in other contexts-- they are examples of what Higman~\cite{higman1988strongly} called ``strongly regular designs'', Neumaier~\cite{neumaier1982regular} called ``$1 \frac{1}{2}$'' designs, and what Bose, Shrikhande, and Singhi~\cite{bose1976edge} called ``partial geometric designs.'' Thus, any construction of a new distance-biregular graph where every vertex has eccentricity four and the girth is four also gives a construction of a new strongly regular/$1 \frac{1}{2}$/partial geometric design.

Delorme~\cite{delorme} gave two infinite families of such graphs, and Van den Akker~\cite{van1990distance} provided another sporadic example.

Godsil and Shawe-Taylor~\cite{distanceRegularised} concluded their seminal paper by saying ``There is a clear need to determine whether the classes of distance-biregular graphs mentioned in this paper exhaust, in any sense, the possibilities.'' In this paper, we address this need by describing the first new constructions of ``uniquely'' distance-biregular graphs in over 30 years.

Section~\ref{preliminaries} introduces the basic definitions and notation of distance-biregular graphs. Section~\ref{examples} collects all the previously known examples of distance-biregular graphs, and Section~\ref{moreProperties} describes some further properties of distance-biregular graphs that we will need going forward, including a new non-existence condition. In Section~\ref{secGenDelorme}, we generalize a construction of Delorme~\cite{delorme} and use this generalization to describe a new distance-biregular graph related to the perp system of Mathon~\cite{DCDHM2002}. In Section~\ref{derivedHyperoval}, we describe a way to derive new distance-biregular graphs as subgraphs of larger distance-biregular graphs, and apply this to obtain a new infinite family of distance-biregular graphs.

\section{Distance-Biregular Graphs}\label{preliminaries}

Distance-biregular graphs have similar algebraic and combinatorial properties to distance-regu\-lar graphs, and it will be convenient to set up common terminology for both.

\subsection{Definitions}
Let $G$ be a graph and let \( \pi = \lsb C_0, C_1, \ldots, C_d \rsb \) be a partition of the vertex set of \( G \). We say that \( \pi \) is \textit{equitable} if for all \( 0 \leq i, j \leq d\),  the number of edges from a vertex in \( C_i \) to vertices in \( C_j \) is independent of the choice of the vertex in \( C_i\).

Let \( \vtxa \in V \lp G \rp \) be a vertex, and let the \textit{eccentricity} \( e \) be the maximum distance from \( \vtxa \) to any other vertex in the graph. For \( 0 \leq i \leq e\),  we let \( N_i \lp \vtxa \rp \) be the set of vertices at distance \( i \) from \( \vtxa\).  Note that \( \lsb N_0(\vtxa),N_1(\vtxa),\ldots,N_e(\vtxa) \rsb \) defines a partition of the vertex set of \( G \), called the \textit{distance partition} of \( \vtxa \). We say that \( \vtxa \) is \textit{locally distance-regular} if the distance partition is equitable.

If \( \vtxb \in N_i \lp \vtxa \rp\), then \( \vtxb \) can only be adjacent to vertices in the cells \( N_{i-1}(\vtxa), N_i(\vtxa) \), and \( N_{i+1} \lp \vtxa \rp\).  In particular, a vertex \( \vtxa \) of eccentricity \( e \) is locally distance-regular if and only if for all \( 0 \leq i \leq e \) and all \( \vtxb \) at distance \( i \) from \( \vtxb\),  the numbers
\[ c_i \lp \vtxa \rp = \abs{N_1(\vtxb) \cap N_{i-1}(\vtxa)}, \]
\[ a_i \lp \vtxa \rp = \abs{N_1(\vtxb) \cap N_{i}(\vtxa)}, \]
\[ b_i \lp \vtxa \rp = \abs{N_1(\vtxb) \cap N_{i+1}(\vtxa)} \]
are well-defined independently of the choice of \( \vtxb\).  We refer to these numbers as the \textit{intersection numbers}.

A graph is \textit{distance-regular} if the distance partition from every vertex is equitable and the intersection numbers \( c_i(\vtxa), a_i(\vtxa) \) and \( b_i \lp \vtxa \rp \) are independent of the choice of \( \vtxa\).  In other words, for every pair of vertices \( \vtxa, \vtxb \) at distance \( i \), we can define a global parameter
\[ c_i := \abs{N_1(\vtxb) \cap N_{i-1}(\vtxa)} \]
that depends only on the distance \( i \) between \( \vtxa \) and \( \vtxb\).  The numbers \( a_i \) and \( b_i \) can be similarly defined independent of the choice of specific vertices \( \vtxa \) and \( \vtxb\).  Distance-regular graphs are a well-studied class of graphs, and more information can be found in the monograph of Brouwer, Cohen, and Neumaier~\cite{bcn} or the more recent survey of Van Dam, Koolen, and Tanaka~\cite{vanDistance}.

When we write a bipartite graph as \( ( Y\cup Z,E) \), it is understood that \( Y \) and \( Z \) are the cells of the partition.
A bipartite graph is \textit{distance-biregular} if the distance partition from every vertex is equitable and the intersection numbers \( c_i \lp \vtxa \rp \) and \( b_i \lp \vtxa \rp \) depend only on which cell of the bipartition \( \vtxa \) lies in. If \( G = \lp \bipartB \cup \bipartC, E \rp \) is a bipartite graph, then \( a_i \lp \vtxa \rp = 0 \) for any vertex \( \vtxa\).  If \( G \) is distance-biregular, we can define global intersection numbers \( c_i^{\bipartB}, c_i^{\bipartC}, b_i^{\bipartB} \), and \( b_i^{\bipartC} \) that only depend on the bipartition.

Locally, distance-regular graphs and distance-biregular graphs behave similarly. However, this local extension of distance-regular graphs to distance-biregular graphs does not extend further, since requiring every vertex to be locally distance-regular is quite restrictive.

\begin{theorem}[Godsil and Shawe-Taylor~\cite{distanceRegularised}]\label{globalLocal}Let \( G \) be a graph where every vertex is locally distance-regular. Then \( G \) is either distance-regular or distance-biregular.
\end{theorem}

\subsection{Generalized Polygons}
A \textit{generalized polygon} is a bipartite graph with girth twice the diameter.  They were introduced by Tits~\cite{tits1959trialite} in studying groups of Lie type and were a particular motivation for Godsil and Shawe-Taylor~\cite{distanceRegularised} because every vertex is locally distance-regular even when the graph is not regular.

\begin{example}\label{genPoly}Let \( G = \lp \bipartB \cup \bipartC, E \rp \) be a bipartite $\lp k, \ell \rp$-semiregular graph with diameter \( d \) and girth \( 2d\),  and let \( \vtxa \in \bipartB \). Let \( \vtxb \) be a vertex at distance \( i \) from \( \vtxa\).  If \( i < d \), then \( \vtxb \) only has one neighbour at distance \( i-1 \) from \( \vtxa\).  Then for all \( 1 \leq i \leq d-1\),  we have \( c_i \lp \vtxa \rp =1\),  and \( b_i \lp \vtxa \rp = k-1 \) if \( i \) is even and \( b_i \lp \vtxa \rp = \ell-1 \) if \( i \) is odd. Similarly, \( c_d = k \) if \( d \) is even and \( c_d = \ell \) if \( d \) is odd, and \( b_0 = k\).

Flipping \( k \) and \( \ell\),  a similar argument holds for \( \vtxc \in \bipartC\),  so \( G \) is distance-biregular.\end{example}

A \text{thick} generalized polygon has minimum degree at least three. Yanushka \cite{yanushka1981order} proved that a generalized polygon that is not thick is the \( k \)-fold subdivision of a multiple edge or the \( k \)-fold subdivision of a thick generalized polygon, and further, thick generalized polygons are semiregular. Feit and Higman~\cite{generalizedPolyThick} proved that any thick generalized polygon has diameter \( d = 2, 3, 4, 6 \), or \( 8\).  Infinite families exist for each of these diameters~\cite{thas1995generalized,van2012generalized}.

A generalized polygon with vertices of degree two can also be distance-biregular if it is the \textit{subdivision graph} obtained by subdividing every edge of a regular generalized polygon exactly once. This leads to the following characterization of Mohar and Shawe-Taylor~\cite{biregularCage} of distance-biregular graph with vertices of valency two.

\begin{theorem}[Mohar and Shawe-Taylor~\cite{biregularCage}]\label{subdivision}Let \( G = \lp \bipartB \cup \bipartC, E \rp \) be a distance-biregular graph where vertices in \( \bipartC \) have valency two. Then \( G \) is either \( K_{2,k} \) or the subdivision graph of a Moore graph or regular generalized polygon.
\end{theorem}

\subsection{Notation and Basic Properties}

Let \( G = \lp \bipartB \cup \bipartC, E \rp \) be a bipartite graph.

Let \( G_2 \) be the graph on the vertex set \( \bipartB \cup \bipartC \) where two vertices are adjacent in \( G_2 \) if and only if they are at distance two in \( G\).  Since \( G \) is bipartite, we see that \( G_2 \) is disconnected. Thus we may let \( H_{\bipartB} \) be the graph \( G_2 \) induced by vertex set \( \bipartB \), and similarly for \( H_{\bipartC}\).  We will refer to \( H_{\bipartB} \) and \( H_{\bipartC} \) as the \textit{halved graphs} induced by \( \bipartB \) and \( \bipartC\),  respectively. Delorme~\cite{delorme} and Mohar and Shawe-Taylor~\cite{biregularCage} observed that the halved graphs of a distance-biregular graph are distance-regular.

The adjacency matrix of a bipartite graph \( G \) has the form
\[ \pmat{\mathbf{0} & N \\ N^T & \mathbf{0}}, \]
where \( N \) is the \( \abs{\bipartB} \times \abs{\bipartC} \) \textit{biadjacency matrix}.

Now suppose that \( G = \lp \bipartB \cup \bipartC, E \rp \) is distance-biregular, with valencies \( k = b_0^{\bipartB} \) and \( \ell = b_0^{\bipartC} \).

One consequence of every vertex in the same cell having the same intersection numbers is that the eccentricity of vertices only depends on the cell of the bipartition they lie in. For \( \bipartS \in \lsb \bipartB, \bipartC \rsb\),  we will denote the maximum eccentricity of vertices in \( \bipartS \) by \( d_{\bipartS}\),  and refer to it as the \textit{covering radius} of \( \bipartS\).  We denote the diameter by \( d\),  and note that at least one of \( d_{\bipartB} \) or \( d_{\bipartC} \) is \( d\).

\begin{lemma}[Delorme~\cite{delorme}]\label{coveringRadii}Let \( G = \lp \bipartB \cup \bipartC, E \rp \) be a distance-biregular graph with diameter \( d \). If \( d_{\bipartB} = d \) then \( d_{\bipartC} \geq d-1\),  and if \( d \) is odd then \( d_{\bipartC} = d\).
\end{lemma}

\begin{lemma}[Delorme~\cite{delorme}]\label{oddRegular}If \( G \) is distance-biregular with odd diameter, then \( G \) is regular.\end{lemma}

The class of regular distance-biregular graphs is equivalent to the class of bipartite distance-regular graphs. In this paper, we are primarily interested in distance-biregular graphs which are not regular, so we will assume the diameter of our graphs is even.

For \( 0 \leq 2i \leq d_{\bipartB} \) we have \( c_{2i}^{\bipartB} + b_{2i}^{\bipartB} = k\),  and for \( 1 \leq 2j+1 \leq d_{\bipartB}\),  we have \( c_{2i+1}^{\bipartB} + b_{2i+1}^{\bipartB} = \ell\).  Following the notation of Delorme~\cite{delorme}, we can compactly express the parameters of a distance-biregular graph in an \textit{intersection array} by
\[ \vmat{k; & c_1^{\bipartB}, & c_2^{\bipartB}, & \ldots, & c_{d_{\bipartB}-1}^{\bipartB} & c_{d_{\bipartB}}^{\bipartB} \\ \ell; & c_1^{\bipartC}, & c_2^{\bipartC}, & \ldots, & c_{d_{\bipartC}-1}^{\bipartC} & c_{d_{\bipartC}}^{\bipartC}} \quad \text{ or } \quad \vmat{k; & c_1^{\bipartB}, & c_2^{\bipartB}, & \ldots, & c_{d_{\bipartB}-1}^{\bipartB} & c_{d_{\bipartB}}^{\bipartB} \\ \ell; & c_1^{\bipartC}, & c_2^{\bipartC}, & \ldots, & c_{d_{\bipartC}}^{\bipartC} &}. \]

Note that \( c_1^{\bipartB} = 1 = c_1^{\bipartC} .\)

\section{Examples}\label{examples}

In this section we describe the known constructions of distance-biregular graphs and relate them to other structures of interest. We may assume without loss of generality that \( 2 \leq k < \ell \), so by Lemma~\ref{oddRegular} we can restrict ourselves to the case where the diameter is even.

\subsection{\texorpdfstring{$d=2$}{complete}}

A complete bipartite graph \( K_{\ell, k} \) is distance-biregular with intersection array
\[ \vmat{k; & 1, & k \\ \ell; & 1, & \ell}. \]
These are clearly the only distance-biregular graphs of diameter two.

\subsection{\texorpdfstring{$d_{\bipartB} = 4, d_{\bipartC} = 3$}{quasi}}

Distance-biregular graphs with \( d_{\bipartB} = 4 \) and \( d_{\bipartC} = 3 \) can be identified with a certain kinds of combinatorial designs.

\begin{theorem}[Delorme~\cite{delorme}, Shawe-Taylor~\cite{shaweTaylor}]\label{quasiDist}A graph \( G \) is distance-biregular with covering radii 4, 3 if and only if it is the incidence graph of a quasi-symmetric design where any two blocks are either disjoint, or they intersect in \( s \) common points.
\end{theorem}

Neumaier~\cite{neumaier1982regular} classified quasi-symmetric designs into four main families, as well as exceptional quasi-symmetric designs. The classes that give rise to distance-biregular graphs are Steiner systems and affine resolvable designs, as well as some of the exceptional quasi-symmetric designs. The known exceptional distance-biregular quasi-symmetric designs are a \( \lp 21, 6, 4 \rp \) quasi-symmetric design coming from the Golay code~\cite{goethals1970strongly} and a \( \lp 22, 6, 5 \rp \) quasi-symmetric design coming from the Witt design~\cite{witt}. More information about quasi-symmetric can be found in the monograph of Shrikhande and Sane~\cite{shrikhande1991quasi}, and a treatment of quasi-symmetric designs through the perspective of distance-biregular graphs can be found in the thesis of Shawe-Taylor~\cite{shaweTaylor}.

\subsection{\texorpdfstring{$d_{\bipartB} = 4, d_{\bipartC} = 4$}{srd}}

A partial geometry \( pg \lp s, t, \alpha \rp \) is a distance-biregular graph with intersection array
\[ \vmat{s+1; & 1, & 1, & \alpha, & s+1 \\ t+1; & 1, & 1, & \alpha, & t+1}. \]

A pg \( \lp s, t, s+1 \rp \) is a Steiner system, and pg \( \lp s, t, 1 \rp \) is a generalized quadrangle. A pg \( \lp s, t, s \rp \) is a \textit{transversal design} of order \( s+1 \)  and degree \( t+1\),  and its existence is equivalent to the existence of \( t-1 \) mutually orthogonal Latin squares of order \( s+1\).  Partial geometries with \( 1 < \alpha < s, t \) are \textit{proper partial geometries}, and there are both sporadic examples and infinite families of proper partial geometries. More information can be found in the surveys of Brouwer and Van Lint~\cite{brouwer1984strongly}, De Clerck and Van Maldeghem~\cite{de1995some}, or Thas~\cite{thas2007partial}.

Quasi-symmetric designs and partial geometries are objects of considerable interest, so the smallest examples of ``uniquely'' distance-biregular graph are distance-biregular graphs with diameter four and girth four. However, genera\-lizations containing these distance-biregular graphs have been introduced under several names as generalizations of other structures of interest.

Bose, Shrikhande, and Singhi~\cite{bose1976edge} coined the term \textit{partial geometric design} as the multigraph analogue to partial geometries.  Bose, Bridges, and Shrikhande \cite{bose1976characterization} further explored the spectral properties of partial geometric designs.

Neumaier~\cite{neumaier1980t12} introduced \( t \frac{1}{2} \)-designs, and showed that the only proper examples came from \( 1 \frac{1}{2} \)-designs. Although the perspective is different, the definition of partial geometric designs and \( 1 \frac{1}{2} \)-designs are equivalent.

Higman~\cite{higman1988strongly} defined a class of coherent configurations as the incidence algebra of what he termed \textit{strongly regular designs}, which he noted contains distance-biregular graphs of diameter four and is contained in the classes of partial geometric or \( 1 \frac{1}{2} \) designs.

\mysqueeze{0.3pt}{We now describe the known constructions from Delorme~\cite{delorme} and Van den Akker \cite{van1990distance} of distance-biregular graphs of diameter four which are not partial geometries.} Let
\[ \qbino{n}{q} = \frac{q^{n}-1}{q-1}. \]

\begin{example}\label{conic}Let \( q \) be a prime power and consider the six-dimensional affine space over \( GF \lp q \rp\).
Consider the (hyperbolic) cone \( X_1 X_2 -X_3X_4 + X_5X_6 = 0 \). Let $\cR^*$ denote the set of three-spaces
in this cone. Then $|\cR^*| = 2(q+1)(q^2+1)$, see~\cite[Lemma 9.4.1]{bcn}.
Pick some $M_0 \in \cR^*$, for instance, $M_0 = \< e_1, e_3, e_5 \>$.
Let \( \cS^* \) denote the set of $M \in \cR^*$ with $\dim(M \cap M_0) \in \{ 1, 3 \}$.
Then $|\cS^*| = (q+1)(q^2+1)$.
Delorme~\cite{delorme} defined a bipartite graph with \( \bipartB \) the \( q^6 \) points of the affine space and \( \bipartC \) the set of three-dimensional affine subspaces parallel to an element of $\cS^*$. Adjacency is incidence.
This graph is distance-biregular with intersection array
  \[ \vmat{\qbino{4}{q}; & 1, & q+1, & q^2, & \qbino{4}{q} \\ q^3; & 1, & q, & q^2 + q, & q^3}. \]
\end{example}

Note that Delorme's description in~\cite{delorme} is slightly incorrect as it uses $\cR^*$ instead of $\cS^*$.
This example was also discussed by Van den Akker~\cite{van1990distance}.

\medskip

Consider a projective plane \( PG \lp 2, q \rp \), and let \( r \) divide \( q \). A \textit{maximal arc} \( \mathcal{A} \) of degree \( r \) is a set of points such that every line meets \( \mathcal{A} \) at 0 or \( r \) points. Denniston~\cite{denniston} and Ball, Blokhuis, and Mazzocca~\cite{ball1998easier, ball1997maximal} proved that a maximal arc exists if and only if \( q \) and \( r \) are both powers of two.

\begin{example}\label{maximalArc}Let \( V \) be the three-dimensional vector space over \( GF \lp q \rp \) and let \( \mathcal{A} \) be a maximal arc in \( PG \lp 2, q \rp \) of degree \( r \). Let \( \hat{\mathcal{A}} \) be the dual of the maximal arc, that is, a set of two-dimensional subspaces of \( V \) such that any one-dimensional subspace of \( V \) is incident with \( 0 \) or \( r \) elements of \( \mathcal{\hat{A}}\).

  Let \( s = |\hat{\mathcal{A}}|\).  We define a bipartite graph where \( \bipartB \) is the \( q^3 \) points of \( V\),  and \( \bipartC \) is the set of the \( qs \) affine cosets of the elements of \( \hat{\mathcal{A}}\). Adjacency is incidence.  Delorme~\cite{delorme} showed this graph is distance-biregular with intersection array
  \[ \vmat{s; & 1, & r, & q \lp s-1 \rp /r, & s \\ q^2; & 1, & q, & s-1, & q^2}. \]
\end{example}

From the work of Denniston~\cite{denniston}, it follows that such a distance-biregular graph exists for \( q, r \) both powers of two. In Section~6.2, Van den Akker~\cite{van1990distance} expanded on Example~\ref{maximalArc} and proved that the existence of such a distance-biregular graph with intersection array
\[ \vmat{n+2; & 1, & 2, & n \lp n-1 \rp /2, & n+2 \\ n^2; & 1, & n, & n+1, & n^2} \]
implies the existence of a projective plane of order \( n\),  which can be used to rule out certain intersection arrays of a similar form to Example~\ref{maximalArc}.

\medskip

The final known example of distance-biregular graphs with diameter four is a sporadic example due to Van den Akker~\cite{van1990distance}.

\begin{example}\label{hallJanko}The Hall-Janko-Wales graph is a strongly regular graph on 100 vertices. It was constructed by Hall and Wales~\cite{hall1968simple}, and key properties are described in Section 10.32 of Brouwer and Van Maldeghem~\cite{brouwer2022strongly}. Let \( \bipartB \) be the vertices of the Hall-Janko-Wales graph and let \( \bipartC \) be the cliques of size 10. Van den Akker~\cite{van1990distance} proved that the bipartite graph with vertex set \( \bipartB \cup \bipartC \) and the incidence relation of inclusion is distance-biregular. This graph has intersection array
  \[ \vmat{28; & 1, & 4, & 6, & 8 \\ 10; & 1, & 2, & 12, & 10}. \]
\end{example}

The other halved graph from this construction was described by Bagchi~\cite{bagchi1992regular}.

\subsection{\texorpdfstring{$d \geq 6$}{largeDiameter}}

There are two known families of distance-biregular graphs with unbounded diameter, which have been described by Delorme~\cite{delorme} and Godsil and Shawe-Taylor~\cite{distanceRegularised}.

\begin{example}\label{biJohnson}Let \( n \geq 2k+2 \) for some \( k \geq 1 \), and let \( S \) be a set of size \( n \). Define a bipartite graph \( G = \lp \bipartB \cup \bipartC, E \rp \) with \( \bipartB \) the subsets of \( S \) of size \( k \) and \( \bipartC \) the subsets of \( S \) of size \( k+1 \), with an edge from \( \vtxa \in \bipartB \) to \( \vtxb \in \bipartC \) if \( \vtxa \subset \vtxb\).  Then \( G \) is distance-biregular with intersection array
  \[ \vmat{n-k; & 1, & 1, & 2, & 2, & \ldots, & k, & k, & k+1 \\ k+1; & 1, & 1, & 2, & 2, & \ldots, & k, & k, & k+1, & k+1}. \]
\end{example}

\begin{example}\label{biGrassmann}Let \( n \geq 2k+2 \) for some \( k \geq 1\).  Let \( q \) be a prime power and let \( V \) be an \( n \)-dimensional vector space over a finite field of \( q \) elements. Define a bipartite graph \( G = \lp \bipartB \cup \bipartC, E \rp \) with \( \bipartB \) the $k$-dimensional subspaces of \( V \) and and \( \bipartC \) the \( \lp k+1 \rp \)-dimensional subspaces of \( V\),  with an edge from \( \vtxa \in \bipartB \) to \( \vtxb \in \bipartC \) if and only if \( \vtxa \) is a subspace of \( \vtxb\).  Then \( G \) is distance-biregular with intersection array
  \[ \vmat{\qbino{n-k}{q}; & \qbino{1}{q}, & \qbino{1}{q}, & \qbino{2}{q}, & \ldots, & \qbino{k}{q}, & \qbino{k}{q}, & \qbino{k+1}{q} \\ \qbino{k+1}{q}; & \qbino{1}{q}, & \qbino{1}{q}, & \qbino{2}{q}, & \ldots, & \qbino{k}{q}, & \qbino{k}{q}, & \qbino{k+1}{q}, & \qbino{k+1}{q}}. \]
\end{example}

The only other known examples of distance-biregular graphs with \( d \geq 6 \) come from generalized hexagons when \( d= 6 \) and generalized octagons when \( d= 8 \). More information can be found in the surveys of Thas~\cite{thas1995generalized} and Van Maldeghem~\cite{van2012generalized}.

\section{Properties}\label{moreProperties}

In this section, we describe some further properties of distance-biregular graphs which will be helpful in constructing new examples.

\subsection{Parameter Relations and Feasibility}

It is well known that many parameters of a distance-regular graph can be determined from the intersection array, and many of these results extend to distance-biregular graphs. One such example is the following result describing the number of vertices in a graph, which we can frame in terms of locally distance-regular vertices.

\begin{lemma}\label{kiS}Let \( \vtxa \) be a locally distance-regular vertex with eccentricity $e$ and let \( k_i \) be the number of vertices at distance \( i \) from \( \vtxa\).  Then \( k_0 = 1 \) and for \( 0 \leq i \leq e-1\),  we have
\[ k_{i+1} = \frac{b_i k_i}{c_{i+1}}. \]
\end{lemma}

This gives us a \textit{feasibility condition} for the intersection array of a distance-biregular graph, since it tells us that for a distance-biregular graph of diameter $d$ to exist, one necessary condition is that for all \( 0 \leq j \leq d-1 \) and \( \bipartS \in \lsb \bipartB, \bipartC \rsb \), the quantity
\[ \prod_{i=0}^{j} \frac{b_i^{\bipartS}}{c_{i+1}^{\bipartS}} \]
must be an integer.

Another feasibility condition comes because the parameters of a distance-biregular graph are not completely independent. It is possible to determine one line of the intersection array from the other by applying the following result due to Delorme~\cite{delorme} or Godsil and Shawe-Taylor~\cite{distanceRegularised}.

\begin{proposition}[Delorme~\cite{delorme}, Godsil and Shawe-Taylor~\cite{distanceRegularised}]\label{delorme2i}Let \( G = \lp \bipartB \cup \bipartC, E \rp \) be a distance-biregular graph of diameter \( d \) with intersection array
  \[ \vmat{k; & 1, & c_2^{\bipartB}, & \ldots, & c_{d_{\bipartB}}^{\bipartB} \\ \ell; & 1, & c_2^{\bipartC}, & \ldots, & c_{d_{\bipartC}}^{\bipartC}}. \]
  For all \( 1 \leq i \leq \lfloor \frac{d-1}{2} \rfloor\),  we have
  \[ c_{2i}^{\bipartB} c_{2i+1}^{\bipartB} = c_{2i}^{\bipartC} c_{2i+1}^{\bipartC} \]
  and
  \[ b_{2i-1}^{\bipartB} b_{2i}^{\bipartB} = b_{2i-1}^{\bipartC} b_{2i}^{\bipartC}. \]
\end{proposition}

\mysqueeze{0.4pt}{Van den Akker~\cite{van1990distance}, Godsil and Shawe-Taylor~\cite{distanceRegularised}, Lato~\cite{phD}, and Secker \cite{secker1989feasibility}} developed sets of feasibility conditions including Proposition~\ref{delorme2i} and Lemma \ref{kiS} to compute tables of feasible parameters. An annotated version of such tables, including the new results from this paper, is given in Section~\ref{paramTable}.

\subsection{The \texorpdfstring{$2$-$\bipartB$}{2-Y}-homogeneous property}

Let \( G = \lp \bipartB \cup \bipartC, E \rp \) be a connected bipartite graph.

Assume that all vertices in \( \bipartB \) have the same covering radius \( d_{\bipartB}  \geq 2\).  Let \( \vtxa \in \bipartB\),  let \( \vtxb \) be at distance two from \( \vtxa\),  and, for \( 1 \leq i \leq d_{\bipartB}\),  let \( \vtxc \) be at distance \( i \) from both \( \vtxa \) and \( \vtxb\).  We introduce a scalar $\gamma_i(\vtxa, \vtxb, \vtxc)$, defined as the count of common neighbours of the vertices $\vtxa$ and $\vtxb$ that are positioned at distance $i-1$ from the vertex $\vtxc$. That is, for \( 1 \leq i \leq d_{\bipartB}\),  we define
\[ \gamma_i(\vtxa, \vtxb, \vtxc)=|N_1(\vtxa)\cap N_1(\vtxb)\cap N_{i-1}(\vtxc)|.  \]

If for some \( 1 \leq i \leq d_{\bipartB} \), the scalar $\gamma_i(\vtxa, \vtxb, \vtxc)$ is invariant for all \( \vtxa \in \bipartB\),  \( \vtxb \) at distance two from \( \vtxa\),  and \( \vtxc \) at distance \( i \) from \( \vtxa \) and \( \vtxb\),  we can define a global constant \( \gamma_i^{\bipartB} := \gamma_i \lp \vtxa, \vtxb, \vtxc \rp\).  If for all \( 1 \leq i \leq d_{\bipartB} -1 \) the quantity \( \gamma_i^{\bipartB} \) is well-defined, we say that \( G \) is \textit{2-$\bipartB$-homogeneous}. If \( \gamma_i^{\bipartB} \) is well-defined for all \( 1 \leq i \leq d_{\bipartB}-2\),  then \( G \) is \textit{almost 2-$\bipartB$-homogeneous}.

To study almost 2-$\bipartB$-homogeneous distance-biregular graphs, Fern\'andez and Penji\'c~\cite{fernandez2023almost} introduced the following scalars.

\begin{definition}\label{delta}
  Let \( G = \lp \bipartB \cup \bipartC, E \rp \) be a distance-biregular graph. For a vertex \( \vtxa \in \bipartB \) and every integer \( 1 \leq i \leq \min \lsb d_{\bipartB}-1, d_{\bipartC}-1 \rsb \), define the scalar \( \Delta_i \lp \bipartB \rp \) by
\[ \Delta_i \lp \bipartB \rp = \begin{cases}\lp b_{i-1}^{\bipartB} -1 \rp \lp c_{i+1}^{\bipartB} -1 \rp - \frac{b_i^{\bipartB} \lp c_{i+1}^{\bipartB}-1 \rp + c_i^{\bipartB} \lp b_{i-1}^{\bipartB}-1 \rp}{c_2^{\bipartB}} \lp c_2^{\bipartC}-1 \rp & i \text{ even} \\ \lp b_{i-1}^{\bipartC} -1 \rp \lp c_{i+1}^{\bipartC} -1 \rp - \frac{b_i^{\bipartC} \lp c_{i+1}^{\bipartC}-1 \rp + c_i^{\bipartC} \lp b_{i-1}^{\bipartC}-1 \rp}{c_2^{\bipartB}} \lp c_2^{\bipartC}-1 \rp & i \text{ odd}
\end{cases}. \]
\end{definition}

\begin{theorem}[Fern\'andez and Penji\'c~\cite{fernandez2023almost}]\label{local}Let \( G = \lp \bipartB \cup \bipartC, E \rp \) be a distance-biregular graph with \( d_{\bipartB} \geq 3 \) and \( \ell \geq 3 \). For \( 2 \leq i \leq \min \{ d_{\bipartB}-1, \allowbreak d_{\bipartC}-1 \}\),  the scalar \( \Delta_i \lp Y \rp = 0 \) if and only if the quantity \( \gamma_i^{\bipartB} \) is well-defined, in which case we have
  \[ \gamma_i^{\bipartB} = \begin{cases} \frac{c_2^{\bipartB} c_i^{\bipartB} \lp b_{i-1}^{\bipartB} -1 \rp}{b_i^{\bipartB} \lp c_{i+1}^{\bipartB} -1 \rp + c_i^{\bipartB} \lp b_{i-1}^{\bipartB}-1 \rp} & i \text{ even} \\ \frac{c_2^{\bipartB} c_i^{\bipartC} \lp b_{i-1}^{\bipartC} -1 \rp}{b_i^{\bipartC} \lp c_{i+1}^{\bipartC} -1 \rp + c_i^{\bipartC} \lp b_{i-1}^{\bipartC}-1 \rp} & i \text{ odd} \end{cases}. \]
\end{theorem}

This gives us a new feasibility conditions, since if \( \Delta_i = 0 \) for some \( i\),  we must have that the above-defined \( \gamma_i \) is a non-negative integer.

\begin{corollary}\label{homFeas}If
  \[ \lp b_{i-1}^{\bipartB} -1 \rp \lp c_{i+1}^{\bipartB} -1 \rp - \frac{b_i^{\bipartB} \lp c_{i+1}^{\bipartB}-1 \rp + c_i^{\bipartB} \lp b_{i-1}^{\bipartB}-1 \rp}{c_2^{\bipartB}} \lp c_2^{\bipartC}-1 \rp = 0, \]
  then
  \[ \frac{c_2^{\bipartB} c_i^{\bipartB} \lp b_{i-1}^{\bipartB} -1 \rp}{b_i^{\bipartB} \lp c_{i+1}^{\bipartB} -1 \rp + c_i^{\bipartB} \lp b_{i-1}^{\bipartB}-1 \rp} \]
  must be a non-negative integer. Similarly, if
  \[ \lp b_{i-1}^{\bipartC} -1 \rp \lp c_{i+1}^{\bipartC} -1 \rp - \frac{b_i^{\bipartC} \lp c_{i+1}^{\bipartC}-1 \rp + c_i^{\bipartC} \lp b_{i-1}^{\bipartC}-1 \rp}{c_2^{\bipartB}} \lp c_2^{\bipartC}-1 \rp = 0, \]
  then
  \[ \frac{c_2^{\bipartB} c_i^{\bipartC} \lp b_{i-1}^{\bipartC} -1 \rp}{b_i^{\bipartC} \lp c_{i+1}^{\bipartC} -1 \rp + c_i^{\bipartC} \lp b_{i-1}^{\bipartC}-1 \rp} \]
  must be a non-negative integer.
\end{corollary}

\subsection{Spectral Excess Theorem}

It is convenient to be able to verify that a graph is distance-biregular without computing the full set of \( 2d+2 \) parameters in the intersection array. Lato~\cite{soloCharacterization} proved one such result with a spectral excess theorem for distance-biregular graphs.

Let \( G \) be a graph with adjacency matrix \( A \), and let \( S \) be a set of vertices. We define an $S$-local inner product of polynomials $f$ and $g$ by
\[ \ip{f}{g}_S = \frac{1}{\abs{S}} \sum_{\vtxa \in S} \ver{\vtxa}^T f \lp A \rp g \lp A \rp \ver{\vtxa}. \]
As noted in Section~2 of~\cite{soloCharacterization}, for a bipartite graph with \( S \) one cell of the bipartition, this inner product is determined by the spectrum and the valencies.

\begin{theorem}[Spectral Excess Theorem~\cite{soloCharacterization}]\label{pBiDistance} Let \( G = \lp \bipartB \cup \bipartC, E \rp \) be a connected \( \lp k, \ell \rp \)-semiregular bipartite graph with diameter \( d \) and \( d+1 \) distinct eigenvalues. Then \( G \) is distance-biregular if and only if there exist orthogonal sequences of polynomials \( p_0^{\bipartB}, \ldots, p_d^{\bipartB} \) and \( p_0^{\bipartC}, \ldots, \allowbreak p_d^{\bipartC} \) such that \( p_d^{\bipartB} \) has degree \( d \) and for every vertex \( \vtxb \in \bipartB \), we have
\[ \norm{p_d^\bipartB}_{\bipartB}^2 =  p_d^{\bipartB} \lp \sqrt{k \ell} \rp = \abs{\lsb \vtxa \in V : d \lp \vtxb, \vtxa \rp = d \rsb}, \]
and \( p_d^{\bipartC} \) has degree \( d \) and for every vertex \( \vtxc \in \bipartC \) we have
\[ \norm{p_d^\bipartC}_{\bipartC}^2 =  p_d^{\bipartC} \lp \sqrt{k \ell} \rp = \abs{\lsb \vtxa \in V : d \lp \vtxc, \vtxa \rp = d \rsb}. \]
\end{theorem}

Applying this result, we can mildly reduce the number of parameters needed to prove that a graph with diameter four is distance-biregular.

\begin{theorem}\label{c3Irrelevant}Let \( G = \lp \bipartB \cup \bipartC, E \rp \) be a bipartite \( \lp k, \ell \rp \)-semiregular graph such that every vertex in \( \bipartB \) is locally distance-regular with parameters
  \[ \vmat{k; & 1, & c_2^{\bipartB}, & c_3^{\bipartB}, & k}. \]
  If \( k > \frac{c_2^{\bipartB} c_3^{\bipartB}}{c_2^{\bipartC}} \) and for every \( \vtxa \in \bipartC \) and \( \vtxb \) at distance two from \( \vtxa \) the number
  \[ c_2^{\bipartC} = |N_1(\vtxa) \cap N_1(\vtxb)| \]
  is independent of the choice of the vertices \( \vtxa \) and \( \vtxb\),  then \( G \) is distance-biregular with intersection array
  \begin{equation}\label{c3Defined} \vmat{k; & 1, & c_2^{\bipartB}, & c_3^{\bipartB}, & k \\ \ell; & 1, & c_2^{\bipartC}, & \frac{c_2^{\bipartB} c_3^{\bipartB}}{c_2^{\bipartC}}, & \ell}.
  \end{equation}
\end{theorem}

\proof Since every vertex in \( \bipartB \) is locally distance-regular, from Section~2 in Lato~\cite{soloCharacterization}, we see that the desired polynomial \( p_4^{\bipartB} \) exists and that vertices in \( \bipartB \) have eigenvalue support of size exactly five, that is, the characteristic vector
of a vertex is a linear combination of exactly five eigenvectors of the
adjacency matrix of \( \bipartB \).

Since the eigenvalues of \( G \) are symmetric, this implies that every vertex in \( \bipartB \) must have zero in its eigenvalue support. It is a well-known result in linear algebra that for any matrix \( N, \) the matrices \( NN^T \) and \( N^TN \) share the same nonzero eigenvalues with multiplicity. In particular, since
\[ A^2 = \pmat{NN^T & \mathbf{0} \\ \mathbf{0} & N^T N}, \]
this implies that \( G \) can only have five distinct eigenvalues.

Now let \( G' = \lp \bipartB' \cup \bipartC', E' \rp \) be a putative distance-biregular graph with the parameters in Equation~\eqref{c3Defined}. For both \( G \) and \( G'\),  every vertex in \( \bipartB \) is locally distance-regular with the same parameters, so from Lemma~\ref{kiS} we have
\[ \abs{\bipartB} = 1 + \frac{k \lp \ell-1 \rp}{c_2^{\bipartB}} + \frac{\lp \ell-1 \rp \lp k-c_2^{\bipartB} \rp \lp \ell-c_3^{\bipartB} \rp}{c_2^{\bipartB} c_3^{\bipartB}} = \abs{\bipartB'} \]
and
\[ \abs{\bipartC} = k + \frac{k \lp \ell-1 \rp \lp k-c_2^{\bipartB} \rp}{c_2^{\bipartB} c_3^{\bipartB}} = \abs{\bipartC'}. \]

We define a sequence of polynomials recursively by
\begin{align*}
  p_0 \lp x \rp &= 1 \\
  p_1 \lp x \rp &= x \\
  p_2 \lp x \rp &= \frac{1}{c_2^{\bipartC}} \lp x p_1 \lp x \rp - \ell \rp \\
  p_3 \lp x \rp &= \frac{c_2^{\bipartC}}{c_2^{\bipartB} c_3^{\bipartB}} \lp x p_2 \lp x \rp - \lp k-1 \rp p_1 \lp x \rp \rp \\
  p_4 \lp x \rp &= \frac{1}{\ell} \lp x p_3 \lp x \rp - \lp \ell-c_2^{\bipartC} \rp p_2 \lp x \rp \rp.
\end{align*}
By construction, these form a sequence of polynomials orthogonal with respect to \( \ip{}{}_{\bipartC'}\),  and further
\[ \ip{p_4}{p_4}_{\bipartC'} = p_4 \lp \sqrt{k \ell} \rp. \] We claim that \( p_4 \) is our desired \( p_4^{\bipartC} \).

As discussed in Section~2.7 of Lato~\cite{phD}, for a bipartite graph, the spectrum can be determined by the spectrum relative to one cell of the bipartition and the sizes of the sets. Further, knowing that every vertex in \( \bipartB \) is locally distance-regular gives us the spectrum relative to \( \bipartB \). Thus \( G \) and \( G' \) are cospectral, and in particular, since they have the same valencies, the inner products \( \ip{}{}_{\bipartC'} \) and \( \ip{}{}_{\bipartC} \) are equivalent. Thus \( p_4 \) does indeed belong to the desired sequence of polynomials orthogonal with respect to \( \ip{}{}_{\bipartC}\),  and further \( \norm{p_4}^{2}_{\bipartC} = p_4 \lp \sqrt{k \ell} \rp\).

We compute that
\[ p_2 \lp \sqrt{k \ell} \rp = \frac{\ell \lp k-1 \rp}{c_2^{\bipartC}} \]
and
\[ p_3 \lp \sqrt{k \ell} \rp = \frac{c_2^{\bipartC}}{c_2^{\bipartB}c_3^{\bipartB}} \lp \sqrt{k \ell} \frac{\ell \lp k-1 \rp}{c_2^{\bipartC}} - \lp k-1 \rp \sqrt{k \ell} \rp =  \frac{\sqrt{k \ell} \lp k-1 \rp \lp \ell - c_2^{\bipartB} \rp}{c_2^{\bipartB} c_3^{\bipartB}}, \]
so
\begin{align*}
 p_4 \lp \sqrt{k \ell} \rp &= \frac{1}{\ell} \lp \frac{k \ell \lp k-1 \rp \lp \ell-c_2^{\bipartC} \rp}{c_2^{\bipartB} c_3^{\bipartB}} - \frac{\lp \ell-c_2^{\bipartC} \rp \ell \lp k-1 \rp}{c_2^{\bipartC}} \rp \\
    &= \frac{\lp k-1 \rp \lp \ell-c_2^{\bipartC} \rp}{c_2^{\bipartB} c_3^{\bipartB}} \lp k - \frac{c_2^{\bipartB} c_3^{\bipartC}}{c_2^{\bipartC}} \rp.
\end{align*}

Now by Lemma~\ref{kiS} we have
\[ \abs{\bipartC} = \abs{\bipartC'} = 1 + \frac{\ell \lp k-1 \rp}{c_2^{\bipartC}} + \frac{\lp k-1 \rp \lp \ell-c_2^{\bipartC} \rp}{c_2^{\bipartB} c_3^{\bipartB}} \lp k-\frac{c_2^{\bipartB} c_3^{\bipartB}}{c_2^{\bipartC}} \rp. \]
Note that for \( G\),  if we fix a vertex \( \vtxa \in \bipartC \), then there are \( \frac{\ell (k-1)}{c_2^{\bipartC}} \) vertices at distance two from \( \vtxa\),  so the number of vertices at distance four from \( \vtxa \) is
\[ \abs{\bipartC} - 1 - \frac{\ell \lp k-1 \rp}{c_2^{\bipartC}} = p_4 \lp \sqrt{k \ell} \rp, \]
and by Theorem~\ref{pBiDistance} we conclude that \( G \) is distance-biregular. By Proposition~\ref{delorme2i} it must have the intersection array in Equation~\eqref{c3Defined}.\qed

\section{Generalization of Delorme's Construction}\label{secGenDelorme}

In this section, we generalize Delorme's construction from Example~\ref{maximalArc} and describe some properties of the generalization.

\subsection{The Generalization}\label{secGeneralization}

A generalization of Delorme's construction is as follows.

\begin{const}\label{constGenDelorme} Let $V$ be a vector space of dimension $n$ over $\mathbb{F}_q$.
 Let $k, s,$ and $d$ be positive integers such that $k \leq n/2$, $s \geq 2$, and $d \geq 2$.
 Let $\cS^*$ be a family of $s$ subspaces of co-dimension $k$ in $V$ with the following properties:
 \begin{enumerate}
  \item for all $v \in V \setminus \{ 0 \}$ we have
  that $|\{M \in \cS^*: v \in M\}| \in \{ 0, d \}$, where both cases occur, and
  \item for all distinct $M, M^* \in \cS^*$ we have $\dim(M \cap M^*) = n-2k$.
 \end{enumerate}
\end{const}

\noindent
We associate Construction \ref{constGenDelorme} with its parameter array $(n, k, q, d, s)$.
We use this to construct a distance-biregular graph generalizing Example~\ref{maximalArc}.

\begin{theorem}\label{thmGenDelorme}Let \( V,k,n,\cS^*, s\),  and \( d \) be as in Construction~\ref{constGenDelorme}. Let \( \bipartB = V\),  let \( \bipartC = \bigcup_{M \in \cS^*} \allowbreak V/M \) and let \( G \) be the bipartite incidence graph on \( \bipartB \cup \bipartC \) with inclusion as the incidence relation. Then \( G \) is distance-biregular with intersection array
\[\vmat{s; & 1, & d, & q^{n-2k}(s-1)/d, & s\\ q^{n-k}; & 1, & q^{n-2k}, & s-1, & q^{n-k}}. \]
\end{theorem}

\proof Every block contains \( q^{n-k} \) points and every point is contained in \( s \) blocks, so \( G \) is \( \lp q^{n-k}, s \rp \)-semiregular.

Let \( B \in \bipartC \) be an arbitrary block. Observe that there exists $M\in\cS^*$ and $b\in V$ such that $B=b+M\in V/M$. If $M'\in\cS^*\setminus\{M\}$, then we have
\[ \dim(M+M')=(n-k)+(n-k)-(n-2k)=n, \]
which implies $M+M'=V$. Thus, if $B'=b'+M'\in V/M'$ with $b'\in V$, then $b-b'=v+v'$
for some $v\in M$ and $v'\in M'$. This implies
\[ b-v=b'+v' \in B\cap B', \]
so $B\cap B'\neq \emptyset$, or equivalently, \( B' \) is at distance two from \( B \). It is clear that $N_2(B)\cap(V/M)=\emptyset$, hence the blocks at distance two from \( B \) are given by
\begin{equation}\label{p27b5}
\bipartC \setminus(V/M) =\bigcup_{M'\in\cS^*\setminus\{M\}} V/M'.
\end{equation}
In particular, if \( B' \) is at distance two from \( B\),  then there exists $M'\in\cS^*\setminus\{M\}$ such that $B'\in V/M'$. Thus the number of points in \( B \cap B' \) is equal to
\[ |M\cap M'| = q^{n-2k}, \]
so \( c_2^{\bipartC} = q^{n-2k}\).

Consider a point \( p \notin B\).  Using Equation~\eqref{p27b5}, the number of blocks at distance two from \( B \) that contain \( p \) is
\[ \sum_{M'\in\cS^*\setminus\{M\}}  |\{B'\in V/M': B'\ni p\}| = |\cS^*\setminus\{M\}| = s-1, \]
so \( c_3^{\bipartC} = s-1\).

Now consider an arbitrary point \( p \in \bipartB\).

If \( p' \) is some other point, the number of blocks that contain both \( p \) and \( p' \) is
\[ \sum_{M'\in\cS^*}|\{B'\in V/M': B'\ni p,p'\}| =  |\{M'\in\cS^*: M'\ni p-p'\}| \in\{0,d\}. \]
It follows that if \( p' \) is at distance two from \( p\),  then there are \( d \) blocks containing both \( p \) and \( p' \). Thus \( c_2^{\bipartB} = d\).

By Theorem~\ref{c3Irrelevant}, we conclude that \( G \) is a distance-biregular graph with the given parameters.\qed

We observe that the resulting distance-biregular graphs do indeed generalize Example~\ref{maximalArc} of Delorme~\cite{delorme}.

\begin{example}\label{exa:maximal_arc}
With reference to Construction~\ref{constGenDelorme}, let $k=1$ and $n=3$. Then $\cS^*$ is a maximal arc, and we
have $(n,k,q,d,s)=(3,1,2^m, 2^{m'}, 2^{m+m'}-2^m+2^{m'})$.
\end{example}

A \textit{projective \( \lp N, K, h_1, h_2 \rp \) set} \( \mathcal{O} \) is a proper, non-empty set of $N$ points of the projective space \( PG \lp K-1, q \rp \) with the property that every hyperplane meets \( \mathcal{O} \) in \( h_1 \) or \( h_2 \) points. Calderbank and Kantor~\cite{MR818812} gave a survey connecting projective sets, two-weight codes, and strongly regular graphs, including tables of known examples.
We give a direct proof connecting projective sets to Construction~\ref{constGenDelorme}.

\begin{theorem}\label{projectiveDel}Let \( V,k,n,\cS^*, s\),  and \( d \) be as in Construction~\ref{constGenDelorme}.  Let \( Y \) be the collection of one-dimensional subspaces \( U \) such that
  \[ | \lsb M \in \cS^* : U \subset M \rsb | = d. \]
  Then \( Y \) is a projective
  \[ \lp \frac{s}{d} \qbino{n-k}{q}, n, \frac{1}{d} \qbino{n-k}{q} + \frac{s-1}{d} \qbino{n-k-1}{q}, \frac{s}{d} \qbino{n-k-1}{q} \rp \]
  set.
\end{theorem}

\proof By the assumption on \( \cS^* \), we know that for \( M_1, M_2 \in \cS^* \), we have \( V = M_1 + M_2\).  Then if \( H \) is a hyperplane, it can contain at most one element of \( \cS^* \). If \( M \in \cS^* \) is not contained in \( H\),  there are \( \qbino{n-k-1}{q} \) points of \( Y \) that lie in a one-dimensional subspace of \( H \). We thus count the pairs \( \lp M, x \rp \) with \( M \in \cS^* \) and \( x \in Y \cap \qbino{H}{q} \) such that \( x \in M \) in two ways to get
\[ d \cdot |Y \cap \qbino{H}{q}| = \begin{cases}\abs{\cS^*} \qbino{n-k-1}{q} & Y \cap \qbin{H}{n-k} = \emptyset \\ \qbino{n-k}{q} + \lp \abs{\cS^*}-1 \rp \qbino{n-k-1}{q} & \text{ otherwise}.\end{cases} \]
Dividing both sides by \( d \) gives the desired result.\qed

It is well known that a two-intersection set gives rise to a strongly regular graph,
for instance, see \cite{MR818812}.
A strongly regular graph $(v, k, \lambda, \mu)$ has eigenvalues $k \geq \tilde{r} \geq 0 > \tilde{s}$ with multiplicities $1$, $f_1$, and $f_2$, respectively.

\begin{corollary}\label{srgProj}Let \( V,k,n,\cS^*, s\),  and \( d \) be as in Construction~\ref{constGenDelorme}, and let \( H_{\bipartB} \) be the graph on vertex set \( V \) with two vertices \( x, y \) adjacent if there exists some \( M \in \cS^* \) such that \( x \) and \( y \) lie in \( V / M \). Then \( H_{\bipartB} \) is strongly regular with parameters
\begin{align*}
  & v = q^n, && \tilde{r} = -\frac{s}{d} + \frac{q^{n-k}}{d},\\
 & k = \frac{s}{d} (q^{n-k}-1), && \tilde{s} = -\frac{s}{d},\\
 & \lambda = \mu -\frac{2s}{d} + \frac{q^{n-k}}{d}, && f_1 = (q-1)h_1,\\
 & \mu = q^{n-2k} \frac{s(s-1)}{d^2,} && f_2 = (q-1)h_2.
\end{align*}
\end{corollary}

Theorem~\ref{thmGenDelorme} can be applied to obtain previously-known distance-regular graphs.

\begin{example}\label{exa:spread}
  With reference to Construction~\ref{constGenDelorme}, suppose $d=1$. Then $\cS^*$ is the dual of a partial spread. As $k \leq n/2$, this forces $n=2k$.
The halved graph \( H_{\bipartB} \) induced by $V$ is a strongly regular graph with Latin square parameters.
In particular, setting $n=2$ and $k=1$, $\cS^*$ may be regarded as a set of directions
of lines in the affine plane $V$. Observe that \( H_{\bipartB} \) is imprimitive if and only if
$\abs{\cS^*}=q^{k}$, that is, one less than the number of $k$-spaces in a spread. In this case we find parameters $(n, k, q, d, s) = (2k, k, q, 1, q^k)$, so \( G \) is distance-regular.
\end{example}

More notably, we can apply Theorem~\ref{thmGenDelorme} to obtain a new distance-biregular graph.

\begin{example}\label{exa:Mathon}
Mathon et al.~\cite{DCDHM2002} describe the dual of a family of $21$
$4$-spaces in a vector space of dimension $6$ over $\mathbb{F}_3$
such that (1) each nonzero vector lies in $0$ or $3$ elements
of the family and (2) the meet of distinct elements of the family is a $2$-space.
In other words, it is an example for Construction \ref{constGenDelorme} with
$(n,k,q,d,s)=(6,2,3,3,21)$.
Bamberg and De Clerck~\cite{BambergDeClerck2010} gave a geometric description.

Using Theorem~\ref{thmGenDelorme}, we obtain a distance-biregular graph with parameters
\[ \vmat{21; & 1, & 3, & 60, & 21 \\ 81; & 1, & 9, & 20, & 81}. \]
\end{example}

\subsection{The Dual Formulation}\label{sec:dual}

The dual of Construction \ref{constGenDelorme} reads as follows:

\begin{const}\label{constGenDelormeDual}
 Let $V$ be a vector space of dimension $n$ over $\mathbb{F}_q$.
 Let $k, s,$ and $d$ be a positive integers such that $k \leq n/2$,
 $s \geq 2$, and $d \geq 2$,
 Let $\cS$ be a family of $s$ subspaces of dimension $k$ in $V$
 with the following properties:
 \begin{enumerate}
  \item for any hyperplane $H$ we have
  that $|\{M \in \cS: M \subseteq H\}| \in \{ 0, d \}$, where both cases occur, and
  \item for all distinct $M, M^* \in \cS$ we have $\dim(M \cap M^*) = 0$.
 \end{enumerate}
\end{const}

Some statements are easier to show in this dual formulation which corresponds to the \textit{dual hyperoval}, see, for example, Yoshiara~\cite{MR1703601}.
A $2$-dimensional dual hyperoval in $PG(4,4)$ gives a family $\cS$ of $3$-spaces of
the $5$-dimensional vector space over $\mathbb{F}_4$
consisting of $22$ members with the desired property.
For a survey, see Dempwolff~\cite{MR4540717} and Yoshiara~\cite{Yoshiara_Pingree}.

A $2$-dimensional dual hyperoval in $PG(4,4)$ is a collection of
$4^2+4+2=22$ planes $\cS$ such that
\begin{itemize}
\item[(a)] Any three distinct members of $\cS$ intersect trivially.
\item[(b)] Any two distinct members of $\cS$ intersect at a projective point.
\end{itemize}
Passing to the dual, in the vector space language, we will have a family $\cS$ of $2$-spaces
which are pairwise trivially intersecting,
and every hyperplane contains at most $2$ members of $\cS$.
In fact, according to Yoshiara \cite{MR1703601}, every hyperplane contains
exactly $0$ or $2$ members of $\cS$. Thus, the objects we are looking for
in the case of $(n,k,d)=(2k+1,k,2)$ are precisely $k$-dimensional dual hyperovals
in $PG(2k,q)$.

\subsection{Restrictions on Parameters}

We show that \( s \) is not an independent parameter.

\begin{lemma}\label{lemCount3a}
 Let \( V,k,n,\cS, s\),  and \( d \) be as in Construction~\ref{constGenDelormeDual}.
 We have
 \begin{align}
  s = \frac{(d-1)(q^{n-k}-1)}{q^{n-2k}-1}+1 = \frac{d(q^{n-k}-1)-q^{n-2k}(q^k-1)}{q^{n-2k}-1}.\label{count3a}
 \end{align}
\end{lemma}
\proof
Counting the pairs \( \lp U, M \rp \) where \( U \) is a one-dimensional subspace of \( V \) and \( M \) is an element of \( \cS \) that contains \( U \), we see
\begin{equation}\label{count1}
  \abs{Y} d = s \qbino{n-k}{q}.
\end{equation}
If we count triples \( \lp U, M_1, M_2 \rp \) where \( M_1 \) and \( M_2 \) are distinct elements in \( \cS^* \) and \( U \) is a one-dimensional subspace of \( V \) in \( M_1 \cap M_2\),  we have
\begin{equation}\label{count2}
\abs{Y} d \lp d-1 \rp = s \lp s-1 \rp \qbino{n-2k}{q}.
\end{equation}

Combining Equation~\eqref{count1} with Equation~\eqref{count2}, we see
\[ \qbino{n-k}{q} \lp d-1 \rp = \lp s-1 \rp \qbino{n-2k}{q}. \]
This concludes the proof. \qed

We can also show that $d$ is necessarily a power of $p$.

\begin{proposition}\label{prop:upperbnd_d}
 Let $V, k, n, \cS, s$, and $d$ be as in Construction~\ref{constGenDelormeDual}.
 Write $q=p^t$, $p$ prime.
 Then $d = q^{n-2k} \cdot p^{-i}$ for some nonnegative integer $i$.
\end{proposition}
\proof Since rational algebraic integers are integers, $-\frac{s}{d}$ and $-\frac{s}{d} + \frac{q^{n-k}}d$ are integers.
 Hence, $d$ is a power of $p$.
 So $d = q^{n-2k} p^{-i}$ for some integer $i$.
 It remains to show that $i$ is nonnegative.

 Recall that $s = 1 + (d-1) \frac{q^{n-k}-1}{q^{n-2k}-1}$.
 All elements of $\cS$ are pairwise disjoint.
 Each element of $\cS$ contains $q^k-1$ nonzero vectors,
 while there are $q^n-1$ nonzero vectors in total.
 Thus,
 \[
  (q^k-1) \left( 1 + (d-1) \frac{q^{n-k}-1}{q^{n-2k}-1} \right) \leq q^n-1.
 \]
 Rearranging for $d-1$ yields
 \[
  d-1 \leq \frac{(q^n-1)(q^{n-2k}-1)}{(q^k-1)(q^{n-k}-1)} - \frac{q^{n-2k}-1}{q^{n-k}-1} < \frac{(q^n-1)(q^{n-2k}-1)}{(q^k-1)(q^{n-k}-1)}.
 \]
 Suppose that $d \geq q^{n-2k} p$. Then
 \[
  p < \frac{(1-q^{-n})(1-q^{-n+2k})}{(1-q^{-k})(1-q^{-n+k})} + q^{-n+2k}.
 \]
 For $p \geq 3$, using $n \geq 2k$, it is easily verified that the right-hand side is less than $3$, a contradiction.
 For $p=2$, we distinguish several cases.
 For $n=2k$, the inequality becomes $2 < 0 + 1$, a contradiction.
 For $n=2k+1$, $1 - q^{-n+2k} < 1 - q^{-k}$ and $1 - q^{-n+k} \geq \frac34$, so the inequality becomes
 \[
  2 < \frac43 (1-q^{-n}) + \frac12 < 2,
 \]
 a contradiction.
 For $n \geq 2k+2$ and either $k \geq 2$ or $q \geq 4$,
 $1 - q^{-k} \geq \frac34$ and $1 - q^{-n+k} \geq \frac78$, so the inequality becomes
 \[
  2 < \frac43 \cdot \frac87 (1-q^{-n}) (1 - q^{-n+k}) + \frac14 < \frac{4 \cdot 8}{3 \cdot 7} + \frac14 < 2,
 \]
 a contradiction.
 It remains the case that $q=2$ and $k=1$. Then the inequality reads
 \[
  2 < \frac{2(2^n-3)}{2^n-2},
 \]
 a contradiction. \qed

\subsection{Bounding \texorpdfstring{$n$}{n}} \label{sec:boundN}

Here we will show that $n$ is bounded by $k$.
For $k=1$, the following proposition is essentially the same as
Lemma 2.4 in \cite{Yoshiara_Pingree}.

\begin{proposition}\label{prop:bound}
Let $V, k, n, \cS, s$, and $d$ be as in Construction~\ref{constGenDelormeDual}.
Suppose
\begin{equation}\label{assume1}
\cS\subsetneqq\qbinom{V}{n-1}\text{ and }d\neq q^{n-2k}.
\end{equation}
Then $n < 4k$.
\end{proposition}
\proof
Observe
\begin{align}
\frac{q^{n-2k}-1}{q^k-1}+1
&=
\frac{q^{n-3k}(q^k-1)+q^{n-3k}+q^k-2}{q^k-1}
\notag\\&>
q^{n-3k}.
\label{625d1}
\end{align}
Since $s\geq2$, \eqref{count3a} implies
\[\frac{(d-1)(q^k-1)}{q^{n-2k}-1}\geq1.\]
Combining this with \eqref{625d1}, we obtain
\begin{equation}\label{625d2}
d>q^{n-3k}.
\end{equation}

On the other hand, by Proposition~\ref{prop:upperbnd_d}
together with the assumption \eqref{assume1}, we have
$d=q^{n-2k}\cdot p^{-j}$
for some positive integer $j$.
Then
\begin{align*}
(p^{j}-1)(q^k-1)
&\equiv(p^{j}-q^{n-2k})(q^k-1)
\pmod{q^{n-2k}-1}
\\&=
-p^j(d-1)(q^k-1)
\\&\equiv0
\pmod{q^{n-2k}-1}
&&\text{(by \eqref{count3a}).}
\end{align*}
Since $j\geq1$, this implies
$(p^{j}-1)(q^k-1)\geq q^{n-2k}-1$, and hence
$p^j>q^{n-3k}$ by \eqref{625d1}.
This in turn implies $d<q^k$.
In view of \eqref{625d2}, we obtain $n<4k$.
\qed
%
%

\section{Derived Hyperovals}\label{derivedHyperoval}

In this section, we give a new family of distance-biregular graphs. This family can be derived as local graphs of Delorme's construction in Example~\ref{maximalArc}, or through a direct geometric argument.

\subsection{Triple intersection numbers}

We begin by describing a method to derive distance-biregular graphs as local graphs of other distance-biregular graphs.

\begin{figure}[ht]
\begin{center}
\begin{tikzpicture}[thick,scale=0.8]
\draw[cNew,very thick] (0,0) -- (3,0);
\draw[cNew,very thick] (0,3) -- (3,3);
\draw[black] (0,6) -- (3,6);
\draw[black] (0,7) -- (3,7);
\draw[cNew,very thick] (6,0) -- (9,0);
\draw[cNew,very thick] (6,4) -- (9,4);
\draw[black] (6,7) -- (9,7);
\draw[cNew,very thick] (0,0) -- (0,3);
\draw[black] (0,3) -- (0,7);
\draw[cNew,very thick] (3,0) -- (3,3);
\draw[black] (3,3) -- (3,7);
\draw[cNew,very thick] (6,0) -- (6,4);
\draw[black] (6,4) -- (6,7);
\draw[cNew,very thick] (9,0) -- (9,4);
\draw[black] (9,4) -- (9,7);
\draw[black] (3,1.5) -- (6,1.5);
\draw[black] (3,4) -- (6,2);
\draw[black] (3,4.5) -- (6,5);
\draw[black] (3,6.5) -- (6,6);
\node () at (1.5,6.5) {$z$};
\node () at (1.5,4.5) {$N_2(z)$};
\node () at (1.5,1.5) {$N_4(z)$};
\node () at (1.5,-0.5) {$\bipartC$};
\node () at (7.5,5.5) {$N_1(z)$};
\node () at (7.5,2) {$N_3(z)$};
\node () at (7.5,-0.5) {$\bipartB$};
\node (c4) at (4,1) {$c^{\bipartC}_4=k$};
\node (b3) at (5.5,1) {$b^{\bipartC}_3$};
\node (c3) at (5.5,3) {$c^{\bipartC}_3$};
\node (b2) at (3.5,3) {$b^{\bipartC}_2$};
\node (c2) at (3.5,5) {$c^{\bipartC}_2$};
\node (b1) at (5.5,4.5) {$b^{\bipartC}_1$};
\node (c1) at (5.5,6.5) {$1$};
\node (b0) at (3.5,6) {$\ell$};
\end{tikzpicture}
\caption{The subgraph induced by $N_3(z)\cup N_4(z)$.}
\label{fig:g3g4z}
\end{center}
\end{figure}
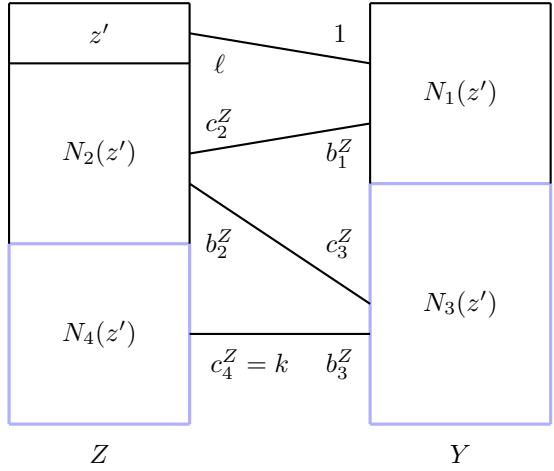

\begin{theorem}\label{tripleIntersect}Let \( G = \lp \bipartB \cup \bipartC, E \rp \) be a distance-biregular graph with intersection array
  \[ \vmat{k; & c_1^{\bipartB}, & c_2^{\bipartB}, & c_3^{\bipartB}, & c_{4}^{\bipartB} \\ \ell; & c_1^{\bipartC}, & c_2^{\bipartC}, & c_3^{\bipartC}, & c_{4}^{\bipartC}}. \]
  In the notation of Definition~\ref{delta}, suppose \( \Delta_3 \lp \bipartB \rp = 0 \) and \[ \gamma_3 = \frac{c_2^{\bipartB}c^{\bipartC}_3(b^{\bipartC}_2-1)}{b^{\bipartC}_3(c^{\bipartC}_4-1)+c^{\bipartC}_3(b^{\bipartC}_2-1)}. \] If \( c_2^{\bipartB} > \gamma_3 \) and \( b_3^{\bipartC} > c_2^{\bipartB}- \gamma_3\),  then for \( z \in \bipartC\),  the subgraph of \( G \) induced by \( N_3 \lp z \rp \cup N_4 \lp z \rp \) is distance-biregular with intersection array
  \begin{equation}\label{tripleArray}
\vmat{b_3^{\bipartC}; & 1, & c_2^{\bipartB}-\gamma_3, & c_3^{\bipartB}, & b^{\bipartC}_3 \\ \ell, & 1, & c^{\bipartC}_2, & \frac{\lp c_2^{\bipartB}-\gamma_3 \rp c_3^{\bipartB}}{c^{\bipartC}_2}, & \ell}.
  \end{equation}
\end{theorem}

\proof Fix \( z \in \bipartC \) and let \( \bipartB' = N_3(z) \), \( \bipartC' = N_4 \lp z \rp \), and \( H \) be the subgraph of \( G \) induced by \( \bipartB' \cup \bipartC' \). It is clear that \( H \) is a \( \lp b_3^{\bipartC}, \ell \rp \)-semiregular graph with diameter at most four. We wish to show that \( H \) is in fact distance-biregular with the parameters detailed in Equation~\eqref{tripleArray}.

Let \( \vtxa \in \bipartB' \) and let \( \vtxb \in \bipartB' \) be at distance two from \( \vtxa\).  Then
\begin{align*}
 \abs{N^{H}_1 \lp \vtxa \rp \cap N^{H}_1 \lp \vtxb \rp} &= \abs{N^{G}_1 \lp \vtxa \rp \cap N^{G}_1 \lp \vtxb \rp \cap N^{G}_4 \lp z \rp} \\
  &= \abs{N^{G}_1 \lp \vtxa \rp \cap N^{G}_1 \lp \vtxb \rp} - \abs{N^{G}_1 \lp \vtxa \rp \cap N^{G}_1 \lp \vtxb \rp \cap N^{G}_2 \lp z \rp} \\& = c_2^{\bipartB}-\gamma_3
\end{align*}
by Theorem~\ref{local}.

Since \( b_3^{\bipartC} > c_2^{\bipartB}-\gamma_3\),  the vertices in \( \bipartB' \) must have eccentricity at least three in \( H\),  so we may choose \( \vtxc \in \bipartC' \) at distance three from \( \vtxa\).  Since every neighbour of \( \vtxc \) in \( G \) is in \( N_3^{G} \lp z \rp\),  we have
\[ \abs{N^{H}_2 \lp \vtxa \rp \cap N^{H}_1 \lp \vtxc \rp} = \abs{N^{G}_2 \lp \vtxa \rp \cap N^{G}_1 \lp \vtxc \rp \cap N^{G}_3 \lp z \rp} = \abs{N^{G}_2 \lp \vtxa \rp \cap N^{G}_1 \lp \vtxc \rp} = c_3^{\bipartB}. \]

Now since \( \ell > c_3^{\bipartB}\),  it follows that every vertex in \( \bipartB' \) is locally distance-regular in \( H \) with eccentricity four. The corresponding line of the intersection array is given by
\[ \vmat{b_3^{\bipartC}; & 1, & c_2^{\bipartB}- \gamma_3, & c_3^{\bipartB}}. \]

Now let \( \vtxa \in \bipartC'\).  Since again every neighbour of \( \vtxa \) in \( G \) is in \( \bipartB'\),  it follows that if \( \vtxb \in \bipartC' \) is at distance two in \( H \) from \( \vtxa\),  then
\[ \abs{N^{H}_1 \lp \vtxa \rp \cap N^{H}_1 \lp \vtxb \rp} = \abs{N^{G}_1 \lp \vtxa \rp \cap N^{G}_1 \lp \vtxb \rp \cap N^{G}_3 \lp z \rp} = \abs{N^{G}_1 \lp \vtxa \rp \cap N^{G}_1 \lp \vtxb \rp} = c_2^{\bipartC}. \]

By Theorem~\ref{c3Irrelevant}, we conclude that \( G \) is distance-biregular with the given parameters.\qed
 \subsection{Another New Construction}

 We can apply Theorem~\ref{tripleIntersect} to get a new family of distance-biregular graphs coming from Example~\ref{maximalArc}.

\begin{theorem}\label{thm:DelormeConstTriple}
  Let \( q = 2^{m} \) for some \( m \geq 2\).  Then there exists a distance-biregular graph with intersection array
  \begin{equation}\label{delormeSub}
\vmat{q+2; & 1, & 2, & \frac{\lp q+1 \rp q}{4}, & q+2 \\ \frac{q \lp q-1 \rp}{2}; & 1, & \frac{q}{2}, & q+1, & \frac{q \lp q-1 \rp}{2}}.
  \end{equation}
\end{theorem}

\proof Example~\ref{maximalArc} with $d=2$ and $q\geq4$ gives the intersection array
\[ \vmat{q^2; & 1, & q, & q+1, & q^2 \\ q+2, & 1, & 2, & \frac{q \lp q+1 \rp}{2}, & q+2}. \]
Then since $b^{\bipartC}_2=\ell-c^{\bipartC}_2=q$ and $b^{\bipartC}_3=k-c^{\bipartC}_3=q(q-1)/2$, we have
\begin{align*}
   \Delta_3(Y) &= (b^{\bipartC}_2-1)(c^{\bipartC}_4-1)-
\frac1{c_2^{\bipartB}}(b^{\bipartC}_3(c^{\bipartC}_4-1)+c^{\bipartC}_3(b^{\bipartC}_2-1))(c^{\bipartC}_2-1)
\\&=
(q-1)((q+2)-1)-
\tfrac{1}{q}(\tfrac{q(q{-}1)}{2}((q{+}2){-}1)+\tfrac{q(q{+}1)}{2}(q{-}1))(2{-}1)
\\&= 0.
\end{align*}
Moreover,
\begin{align*}
\gamma_3&=\frac{c_2^{\bipartB}c^{\bipartC}_3(b^{\bipartC}_2-1)}{b^{\bipartC}_3(c^{\bipartC}_4-1)+c^{\bipartC}_3(b^{\bipartC}_2-1)}
=
\frac{q}{2} <q = c_2^{\bipartB}
\end{align*}
and
\begin{align*}
b^{\bipartC}_3&=k-c^{\bipartC}_3
=
q^2-\frac{q(q+1)}{2}
=
\frac{q(q-1)}{2}
\\&>
\frac{q(q+1)}{4}
=
(c_2^{\bipartB}-\gamma_3)\frac{q+1}{2}
=
c_2^{\bipartB}-\gamma_3,
\end{align*}
and so the conditions of Theorem~\ref{tripleIntersect} are satisfied.\qed

\noindent
We can obtain this same family through a more direct geometric construction.

\begin{const}
\label{constHyperoval}
Let \( q= 2^{m} \) for some \( m \geq 2\).  Fix a point \( x \in \mathbb{F}_q^3 \) and a plane \( \pi \) at infinity. Let \( H^* \) be a set of \( q+2 \) projective lines in \( \pi \) such that every point lies in 0 or 2 lines. We partition \( \pi \) into the points lying in no lines in \( H^* \), called the \textit{exterior points} and the points lying in two lines in \( H^* \), called \textit{interior points}.
Let \( \bipartB \) be the set of \( y  \in \mathbb{F}_q^3 \) such that \( \<x, y\> \cap \pi \) is an exterior point. Let \( \bipartC \) be the set of affine cosets of the elements of \( H^* \) that do not contain \( x \).  Define a bipartite incidence graph \( G \) on vertex set \( \bipartB \cup \bipartC \).\end{const}

\begin{theorem}\label{thmHyperoval}Let \( q = 2^{m} \) for some \( m \geq 2\).  Then Construction~\ref{constHyperoval} is a distance-biregular graph with intersection array \eqref{delormeSub}.
\end{theorem}

\proof If we choose two lines in \( H^*\),  they intersect in a unique interior point, so there are \( \frac{1}{2} \lp q+2 \rp \lp q+1 \rp \) interior points. Thus there are \( \frac{1}{2} q \lp q-1 \rp \) exterior points.
Note that \( \abs{\bipartB} = \frac{1}{2} q \lp q-1 \rp^2 \) and \( \abs{\bipartC} = \lp q+2 \rp \lp q-1 \rp\).

If \( \vtxa \in \bipartB \) and \( L \) is a line in \( H^*\),  then \( \< \vtxa, L \> \) does not contain \( x\),  so \( \vtxa \) is incident to \( q+2 \) blocks in \( \bipartC\).  Now let \( \vtxb \in \bipartC \) and let \( y \in \mathbb{F}_q^3 \setminus \lsb x \rsb \) be contained in \( \vtxb\).  Then \( y \in \bipartB \) precisely when \( \<x,y \> \cap \pi \) is an exterior point, so there are \( \frac{1}{2} q \lp q-1 \rp \) such points \( y \in \bipartB \) incident to \( \vtxb\).  Thus \( G \) is \( \lp q+2, \frac{1}{2} q \lp q-1 \rp \rp \)-semiregular.

Fix $\vtxb \in \bipartC$.
Let \( \vtxc \) be a block that is not parallel to \( \vtxb\),  so \( \vtxb \) and \( \vtxc \) are at distance two.
Then \( \vtxb \) and \( \vtxc \) intersect in an affine line \( L \).
Consider the projection of $L$ from $x$ onto infinity, that is,
$L' = \< x, L \> \cap \pi$.
Each of the other \( q \) projective lines in \( H^* \) must meet \( L' \)
in precisely one point, and each of the points \( y \in L' \) that
lie on some line in \( H^* \) must lie in precisely two projective
lines of \( H^* \). Thus \( \vtxb \) and \( \vtxc \) intersect
in \( \frac{q}{2} \) points in \( \bipartB\).
Hence, $c_2^{\bipartC} = \frac{q}2$.

Now let $\vtxa$ be a point that is not on $\vtxb$.
For $L \in H^*$
we have that $\< \vtxa, L \>$ and $\vtxb$ are at distance $2$
precisely when $L \neq \vtxb \cap \pi$. Hence,
$c_3^{\bipartC} = q+1$.

Now fix $\vtxa \in \bipartB$.
Let $\vtxc \in \bipartB$ be at distance $2$ from $\vtxa$.
Then $\< \vtxa, \vtxc \> \cap \pi$ is an interior point $P$.
Let $L_1, L_2$ denote the two lines of $\pi$ through $P$.
The only common neighbours of $\vtxa$ and $\vtxc$
are $\< \vtxa, L_1 \>$ and $\< \vtxa, L_2\>$.
Hence, $c_2^\bipartB = 2$, so by Theorem~\ref{c3Irrelevant}, the graph is distance-biregular with the specified parameters.\qed

\mysqueeze{0.31pt}{The nontrivial strongly regular halved graph was defined by Huang, Huang, and Lin \cite{huang2007class},} and the set-up used is similar to the construction of Brouwer~\cite{brouwerStronglyHyperoval} and Brouwer, Ihringer, and Kantor~\cite{brouwer2023strongly} to describe the complement.

\section{Feasible Parameters}\label{paramTable}

In the following table, we list all the feasible parameters of non-regular distance-biregular graphs with \( d_{\bipartB} = d_{\bipartC} = 4 \), girth four, and at most 1300 vertices. We also list the parameters of the strongly regular halved graphs.

Using the colour scheme of~\cite{srgDatabase}, green represents that the distance-biregular or strongly regular graph is known to exist, red indicates that it is known not to exist, and yellow represents that the existence is unknown. The new constructions from this paper are marked in blue.

The feasibility conditions used here are those of Section~3.6 of Lato~\cite{phD}, plus the Krein condition which can be found in Delorme~\cite{delorme}.
We then apply Corollary~\ref{homFeas} to establish a new non-existence criterion for distance-biregular graphs. In particular, this criterion rules out the existence of certain distance-biregular graphs, thereby providing a new feasibility condition for this class of graphs.

For more information on the strongly regular halved graphs, the reader is referred to Brouwer and Van Maldeghem~\cite{brouwer2022strongly}. For the distance-biregular graphs, both external and internal references are included.


{\small
\begin{longtable}{| c | c| c |}
  \hline Intersection Array & Halved Graph & Notes \\ \hline
  \endhead
  \cellcolor{cExist}
& \cellcolor{cExist}$(64, 45, 32, 30)$ & Delorme~\cite{delorme} \\ \multirow{-2}{*}{$\vmat{6; & 1, & 2, & 10, & 6 \\ 16; &1, & 4, & 5, & 16}$ \cellcolor{cExist}}  & \cellcolor{cExist} $(24, 20, 16, 20)$ & Ex.~\ref{maximalArc}: $q=4, r=2$ \\ \hline
\cellcolor{cExist}
 & \cellcolor{cExist} $(120, 56, 28, 24)$ & Delorme~\cite{delorme} \\ \multirow{-2}{*}{$\vmat{8; & 1, & 2, & 6, & 8 \\ 15; &1, & 3, & 4, & 15}$} \cellcolor{cExist} & \cellcolor{cExist} $(64, 35, 18, 20)$ & Ex.~\ref{conic}: $q=2$ \\ \hline
\cellcolor{cNew}
 & \cellcolor{cExist} $(196, 135, 94, 90)$ & Constr.~\ref{constHyperoval} \\ \multirow{-2}{*}{$\vmat{10; & 1, & 2, & 18, & 10 \\ 28; &1, & 4, & 9, & 28}$} \cellcolor{cNew} & \cellcolor{cExist} $(70, 63, 56, 63)$ & $q=2$\\ \hline
\cellcolor{cNot}
 & \cellcolor{cUnknown} $(216, 140, 94, 84)$ & Van den Akker~\cite{van1990distance} \\ \multirow{-2}*{$\vmat{8; & 1, & 2, & 21, & 8 \\ 36; &1, & 6, & 7, & 36}$} \cellcolor{cNot} & \cellcolor{cExist} $(48, 42, 36, 42)$ & Section~6.2\\ \hline
\cellcolor{cUnknown}
 & \cellcolor{cExist} $(216, 175, 142, 140)$ & Only known SRG \cite{CRS2018} \\ \multirow{-2}*{$\vmat{15; & 1, & 3, & 28, & 15 \\ 36; &1, & 6, & 14, & 36}$} \cellcolor{cUnknown} & \cellcolor{cExist} $(90, 84, 78, 84)$ & does not work  \\ \hline
\cellcolor{cNot}
 & \cellcolor{cUnknown} $(225, 176, 139, 132)$ & Corollary~\ref{homFeas} \\ \multirow{-2}*{$\vmat{12; & 1, & 3, & 33, & 12 \\ 45; &1, & 9, & 11, & 45}$} \cellcolor{cNot} & \cellcolor{cExist} $(60, 55, 50, 55)$ & $\gamma_2 = \frac{9}{5}$ \\ \hline
\cellcolor{cExist}  & \cellcolor{cExist} $(280, 135, 70, 60)$ & Van den Akker~\cite{van1990distance} \\ \multirow{-2}*{$\vmat{10; & 1, & 2, & 12, & 10 \\ 28; &1, & 4, & 6, & 28}$} \cellcolor{cExist} & \cellcolor{cExist} $(100, 63, 38, 42)$ & Ex.~\ref{hallJanko} \\ \hline
\cellcolor{cUnknown}  & \cellcolor{cUnknown} $(288, 175, 110, 100)$ & \\ \multirow{-2}*{$\vmat{15; & 1, & 3, & 20, & 15 \\ 36; &1, & 6, & 10, & 36}$} \cellcolor{cUnknown} & \cellcolor{cUnknown} $(120, 84, 58, 60)$ & \\ \hline
\cellcolor{cUnknown}  & \cellcolor{cUnknown} $(378, 182, 91, 84)$ & \\ \multirow{-2}*{$\vmat{14; & 1, & 2, & 12, & 14 \\ 27; &1, & 3, & 8, & 27}$} \cellcolor{cUnknown} & \cellcolor{cExist} $(196, 117, 68, 72)$ & \\ \hline
\cellcolor{cUnknown}  & \cellcolor{cUnknown} $(400, 273, 188, 182)$ & \\ \multirow{-2}*{$\vmat{14; & 1, & 2, & 26, & 14 \\ 40; &1, & 4, & 13, & 40}$} \cellcolor{cUnknown} & \cellcolor{cExist} $(140, 130, 120, 130)$ & \\ \hline
\cellcolor{cExist}  & \cellcolor{cExist} $(512, 315, 202, 180)$ & Delorme~\cite{delorme} \\ \multirow{-2}*{$\vmat{10; & 1, & 2, & 36, & 10 \\ 64; &1, & 8, & 9, & 64}$} \cellcolor{cExist} & \cellcolor{cExist} $(80, 72, 64, 72)$ & Ex.~\ref{maximalArc}: $q=4, r=2$ \\ \hline
\cellcolor{cExist}  & \cellcolor{cExist} $(512, 441, 380, 378)$ & Delorme~\cite{delorme} \\ \multirow{-2}*{$\vmat{28; & 1, & 4, & 54, & 28 \\ 64; &1, & 8, & 27, & 64}$} \cellcolor{cExist} & \cellcolor{cExist} $(224, 216, 208, 216)$ & Ex.~\ref{conic}: $q=4$ \\ \hline
\cellcolor{cUnknown}  & \cellcolor{cExist} $(540, 264, 138, 120)$ & \\ \multirow{-2}*{$\vmat{12; & 1, & 2, & 20, & 12 \\ 45; &1, & 5, & 8, & 45}$} \cellcolor{cUnknown} & \cellcolor{cExist} $(144, 99, 66, 72)$ & \\ \hline
\cellcolor{cUnknown}  & \cellcolor{cUnknown} $(560, 273, 140, 126)$ & \\ \multirow{-2}*{$\vmat{14; & 1, & 2, & 18, & 14 \\ 40; &1, & 4, & 9, & 40}$} \cellcolor{cUnknown} & \cellcolor{cExist} $(196, 130, 84, 90)$ & \\ \hline
\cellcolor{cUnknown}  & \cellcolor{cExist} $(560, 351, 222, 216)$ &  \\ \multirow{-2}*{$\vmat{27; & 1, & 3, & 24, & 27 \\ 40; &1, & 4, & 18, & 40}$} \cellcolor{cUnknown} & \cellcolor{cExist} $(378, 260, 178, 180)$ & \\  \hline
\cellcolor{cNot}  & \cellcolor{cUnknown} $(576, 475, 394, 380)$ & Corollary~\ref{homFeas} \\ \multirow{-2}*{$\vmat{20; & 1, & 4, & 76, & 20 \\ 96; &1, & 16, & 19, & 96}$} \cellcolor{cNot} & \cellcolor{cExist} $(120, 114, 108, 114)$ & $\gamma_2 = \frac{8}{3}$ \\ \hline
\cellcolor{cUnknown}  & \cellcolor{cUnknown} $(640, 441, 308, 294)$ & \\ \multirow{-2}*{$\vmat{28; & 1, & 4, & 42, & 28 \\ 64; &1, & 8, & 21, & 64}$} \cellcolor{cUnknown} & \cellcolor{cUnknown} $(280, 216, 166, 168)$ & \\ \hline
\cellcolor{cUnknown}  & \cellcolor{cUnknown} $(676, 459, 314, 306)$ & \\ \multirow{-2}*{$\vmat{18; & 1, & 2, & 34, & 18 \\ 52; &1, & 4, & 17, & 52}$} \cellcolor{cUnknown} & \cellcolor{cExist} $(234, 221, 208, 221)$ & \\ \hline
\cellcolor{cUnknown}  & \cellcolor{cUnknown}$(726, 455, 292, 273)$ & \\* \multirow{-2}*{$\vmat{14; & 1, & 2, & 39, & 14 \\ 66; &1, & 6, & 13, & 66}$} \cellcolor{cUnknown} & \cellcolor{cExist} $(154, 143, 132, 143)$ & \\ \hline
\cellcolor{cUnknown}  & \cellcolor{cUnknown} $(726, 585, 472, 468)$ & \\* \multirow{-2}*{$\vmat{27; & 1, & 3, & 52, & 27 \\ 66; &1, & 6, & 26, & 66}$} \cellcolor{cUnknown} & \cellcolor{cExist} $(297, 286, 275, 286)$ & \\ \hline
\cellcolor{cNew}  & \cellcolor{cExist} $(729, 560, 433, 420)$ & Constr.~\ref{constGenDelorme} $q=3$ \\ \multirow{-2}*{$\vmat{21; & 1, & 3, & 60, & 21 \\ 81; &1, & 9, & 20, & 81}$} \cellcolor{cNew} & \cellcolor{cExist} $(189, 180, 171, 180)$ & $n=6, k=2, d=3, s =21$ \\ \hline
\cellcolor{cUnknown}  & \cellcolor{cExist} $(729, 182, 55, 42)$ & \\ \multirow{-2}*{$\vmat{14; & 1, & 2, & 6, & 14 \\ 27; &1, & 3, & 4, & 27}$} \cellcolor{cUnknown} & \cellcolor{cExist} $(378, 117, 36, 36)$ & \\ \hline
\cellcolor{cUnknown}  & \cellcolor{cExist} $(780, 380, 190, 180)$ & \\ \multirow{-2}*{$\vmat{20; & 1, & 2, & 18, & 20 \\ 39; &1, & 3, & 12, & 39}$} \cellcolor{cUnknown} & \cellcolor{cUnknown} $(400, 247, 150, 156)$ & \\ \hline
\cellcolor{cUnknown}  & \cellcolor{cUnknown} $(875, 552, 355, 336)$ & \\ \multirow{-2}*{$\vmat{24; & 1, & 3, & 42, & 24 \\ 70; &1, & 7, & 18, & 70}$} \cellcolor{cUnknown} & \cellcolor{cUnknown} $(300, 230, 175, 180)$ & \\ \hline
\cellcolor{cUnknown}  & \cellcolor{cUnknown} $(924, 455, 238, 210)$ & \\ \multirow{-2}*{$\vmat{14; & 1, & 2, & 30, & 14 \\ 66; &1, & 6, & 10, & 66}$} \cellcolor{cUnknown} & \cellcolor{cExist} $(196, 143, 102, 110)$ & \\ \hline
\cellcolor{cUnknown}  & \cellcolor{cUnknown} $(924, 650, 460, 450)$ & \\ \multirow{-2}*{$\vmat{40; & 1, & 4, & 45, & 40 \\ 66; &1, & 6, & 30, & 66}$} \cellcolor{cUnknown} & \cellcolor{cUnknown} $(560, 429, 328, 330)$ & \\ \hline
\cellcolor{cUnknown}  & \cellcolor{cUnknown} $(936, 459, 234, 216)$ & \\ \multirow{-2}*{$\vmat{18; & 1, & 2, & 24, & 18 \\ 52; &1, & 4, & 12, & 52}$} \cellcolor{cUnknown} & \cellcolor{cExist} $(324, 221, 148, 156)$ & \\ \hline
\cellcolor{cUnknown}  & \cellcolor{cUnknown} $(945, 560, 343, 315)$ & \\ \multirow{-2}*{$\vmat{21; & 1, & 3, & 45, & 21 \\ 81; &1, & 9, & 15, & 81}$} \cellcolor{cUnknown} & \cellcolor{cUnknown} $(245, 180, 131, 135)$ &\\ \hline
\cellcolor{cNot}  & \cellcolor{cUnknown} $(960, 714, 538, 510)$ & Corollary~\ref{homFeas} \\ \multirow{-2}*{$\vmat{18; & 1, & 3, & 85, & 18 \\ 120; &1, & 15, & 17, & 120}$} \cellcolor{cNot} & \cellcolor{cExist} $(144, 136, 128, 136)$ & $\gamma_2 = \frac{15}{8}$ \\ \hline
\cellcolor{cUnknown}  & \cellcolor{cUnknown} $(960, 833, 724, 714)$ & \\ \multirow{-2}*{$\vmat{35; & 1, & 5, & 102, & 35 \\ 120; &1, & 15, & 34, & 120}$} \cellcolor{cUnknown} & \cellcolor{cExist} $(280, 272, 264, 272)$ & \\ \hline
\cellcolor{cNot}  & \cellcolor{cUnknown} $(1000, 594, 368, 330)$ & Van den Akker~\cite{van1990distance} \\ \multirow{-2}*{$\vmat{12; & 1, & 2, & 55, & 12 \\ 100; &1, & 10, & 11, & 100}$} \cellcolor{cNot} & \cellcolor{cExist} $(120, 110, 100, 110)$ & Section~6.2 \\ \hline
\cellcolor{cUnknown}  & \cellcolor{cUnknown} $(1000, 891, 794, 792)$ & \\ \multirow{-2}*{$\vmat{45; & 1, & 5, & 88, & 45 \\ 100; &1, & 10, & 44, & 100}$}\cellcolor{cUnknown} & \cellcolor{cExist} $(450, 440, 430, 440)$ & \\ \hline
\cellcolor{cUnknown}  & \cellcolor{cExist} $(1024, 693, 472, 462)$ & Constr.~\ref{constGenDelorme} $q \in \{ 2, 4\}$? \\ \multirow{-2}*{$\vmat{22; & 1, & 2, & 42, & 22 \\ 64; &1, & 4, & 21, & 64}$} \cellcolor{cUnknown} & \cellcolor{cExist} $(352, 336, 320, 336)$ & If ex.~$q=4$: $|\mathrm{Aut}| = 2^m$. \\  \hline
\cellcolor{cUnknown}  & \cellcolor{cUnknown} $(1056, 650, 406, 390)$ & \\ \multirow{-2}*{$\vmat{40; & 1, & 4, & 39, & 40 \\ 66; &1, & 6, & 26, & 66}$}\cellcolor{cUnknown} & \cellcolor{cUnknown} $(640, 429, 288, 286)$ & \\ \hline
\cellcolor{cExist}  & \cellcolor{cExist} $(1080, 351, 126, 108)$ & Delorme~\cite{delorme}\\ \multirow{-2}*{$\vmat{27; & 1, & 3, & 12, & 27 \\ 40; &1, & 4, & 9, & 40}$} \cellcolor{cExist} & \cellcolor{cExist} $(729, 260, 97, 90)$ & Ex.~\ref{conic}: $q=3$\\ \hline
\cellcolor{cUnknown}  & \cellcolor{cUnknown} $(1200, 891, 666, 648)$ & \\ \multirow{-2}*{$\vmat{45; & 1, & 5, & 72, & 45 \\ 100; &1, & 10, & 36, & 100}$} \cellcolor{cUnknown} & \cellcolor{cUnknown} $(540, 440, 358, 360)$ & \\ \hline
\cellcolor{cUnknown}  & \cellcolor{cUnknown} $(1210, 585, 296, 270)$ & \\ \multirow{-2}*{$\vmat{27; & 1, & 3, & 30, & 27 \\ 66; &1, & 6, & 15, & 66}$}\cellcolor{cUnknown} & \cellcolor{cUnknown} $(495, 286, 165, 165)$ & \\ \hline
\cellcolor{cNot}  & \cellcolor{cUnknown} $(1225, 1044, 893, 870)$ & Corollary~\ref{homFeas} \\* \multirow{-2}*{$\vmat{30; & 1, & 5, & 145, & 30 \\ 175; &1, & 25, & 29, & 175}$} \cellcolor{cNot} & \cellcolor{cExist} $(210, 203, 196, 203)$ & $\gamma_2 = \frac{25}{7}$ \\ \hline
\end{longtable}}

\section*{Acknowledgements}
We thank the reviewers for their valuable comments.


%
%

\bibliographystyle{alphaurl}
\bibliography{references}

\begin{thebibliography}{DCDHM02}


\bibitem[A90]{van1990distance}
J.~van~den Akker.
\newblock Distance-biregular graphs.
\newblock Master's thesis, Eindhoven University of Technology, 1990.

\bibitem[Bag92]{bagchi1992regular}
B.~Bagchi.
\newblock A regular two-graph admitting the {Hall-Janko-Wales} group.
\newblock {\em Sankhy{\=a}: The Indian Journal of Statistics, Series B}, pages
  35--45, 1992.

\bibitem[BB98]{ball1998easier}
S.~Ball and A.~Blokhuis.
\newblock An easier proof of the maximal arcs conjecture.
\newblock {\em Proceedings of the American Mathematical Society},
  126(11):3377--3380, 1998.

\bibitem[BBM97]{ball1997maximal}
S.~Ball, A.~Blokhuis, and F.~Mazzocca.
\newblock Maximal arcs in {Desarguesian} planes of odd order do not exist.
\newblock {\em Combinatorica}, 17(1):31--41, 1997.

\bibitem[BBS76]{bose1976characterization}
R.~Bose, W.~G. Bridges, and M.~S. Shrikhande.
\newblock A characterization of partial geometric designs.
\newblock {\em Discrete Mathematics}, 16(1):1--7, 1976.

\bibitem[BCN89]{bcn}
A.~E. Brouwer, A.~M. Cohen, and A.~Neumaier.
\newblock {\em Distance-regular graphs}.
\newblock Springer Berlin Heidelberg, 1989.

\bibitem[BDC10]{BambergDeClerck2010}
J.~Bamberg and F.~De~Clerck.
\newblock A geometric construction of {Mathon's} perp-system from four lines of
  {PG}(5, 3).
\newblock {\em Journal of Combinatorial Designs}, 18(6):450--461, 2010.

\bibitem[BIK23]{brouwer2023strongly}
A.~E. Brouwer, F.~Ihringer, and W.~M. Kantor.
\newblock Strongly regular graphs satisfying the 4-vertex condition.
\newblock {\em Combinatorica}, 43(2):257--276, 2023.

\bibitem[Bos63]{bose1963strongly}
R.~Bose.
\newblock Strongly regular graphs, partial geometries and partially balanced
  designs.
\newblock {\em Pacific Journal of Mathematics}, 13(2):389--419, 1963.

\bibitem[Broa]{srgDatabase}
A.~E. Brouwer.
\newblock Parameters of strongly regular graphs.
\newblock \url{https://aeb.win.tue.nl/graphs/srg/srgtab.html}.
\newblock Accessed: 16 April 2025.

\bibitem[Brob]{brouwerStronglyHyperoval}
A.~E. Brouwer.
\newblock Strongly regular graphs from hyperovals.
\newblock \url{https://www.win.tue.nl/~aeb/preprints/hhl.pdf}.
\newblock Accessed: 01 May 2025.

\bibitem[BSS76]{bose1976edge}
R.~Bose, S.~Shrikhande, and N.~Singhi.
\newblock Edge regular multigraphs and partial geometric designs with an
  application to the embedding of quasi-residual designs.
\newblock {\em Colloquie Internazionale sulle Teorie Combinatorie}, 1:49--81,
  1976.

\bibitem[BvL84]{brouwer1984strongly}
A.~E. Brouwer and J.~H. van Lint.
\newblock Strongly regular graphs and partial geometries.
\newblock In {\em Enumeration and Design 85}, page 122. 1984.

\bibitem[BVM22]{brouwer2022strongly}
A.~E. Brouwer and H.~Van~Maldeghem.
\newblock {\em Strongly regular graphs}, volume 182.
\newblock Cambridge University Press, 2022.

\bibitem[CDS80]{cvetkovic}
D.~Cvetkovi{\'c}, M.~Doob, and H.~Sachs.
\newblock {\em Spectra of {Graphs}: Theory and Applications}.
\newblock Academic Press, 1980.

\bibitem[CK86]{MR818812}
R.~Calderbank and W.~M. Kantor.
\newblock The geometry of two-weight codes.
\newblock {\em Bulletin of the London Mathematical Society}, 18(2):97--122,
  1986.

\bibitem[CR{\v S}18]{CRS2018}
D.~Crnkovi{\'c}, S.~Rukavina, and A.~{\v S}vob.
\newblock New strongly regular graphs from orthogonal groups {$O^+(6,2)$} and
  {$O^-(6,2)$}.
\newblock {\em Discrete Mathematics}, 341(10):2723--2728, 2018.

\bibitem[DKT16]{vanDistance}
E.~R. van Dam, J.~Koolen, and H.~Tanaka.
\newblock Distance-regular graphs.
\newblock {\em The Electronic Journal of Combinatorics}, page DS22, 2016.

\bibitem[DCDHM02]{DCDHM2002}
F.~De~Clerck, M.~Delanote, N.~Hamilton, and R.~Mathon.
\newblock Perp-systems and partial geometries.
\newblock {\em Advances in Geometry}, 2(1):1--12, 2002.

\bibitem[DCVM95]{de1995some}
F.~De~Clerck and H.~Van~Maldeghem.
\newblock Some classes of rank 2 geometries.
\newblock In {\em Handbook of Incidence Geometry}, pages 433--473. Elsevier,
  1995.

\bibitem[Del83]{delormeFrench}
C.~Delorme.
\newblock Regularit{\'e} m{\'e}trique forte.
\newblock Technical Report 156, Universit{\'e} Paris-Sud, Orsay, 1983.

\bibitem[Del94]{delorme}
C.~Delorme.
\newblock Distance biregular bipartite graphs.
\newblock {\em European Journal of Combinatorics}, 15(3):223--238, 1994.

\bibitem[Dem19]{MR4540717}
U.~Dempwolff.
\newblock Dimensional dual hyperovals---an updated survey.
\newblock In {\em International Conference on Algebra and Related Topics with
  Applications}, pages 115--142. Springer, 2019.

\bibitem[Den69]{denniston}
R.~Denniston.
\newblock Some maximal arcs in finite projective planes.
\newblock {\em Journal of Combinatorial Theory}, 6(3):317--319, 1969.

\bibitem[FH64]{generalizedPolyThick}
W.~Feit and D.~G. Higman.
\newblock The nonexistence of certain generalized polygons.
\newblock {\em Journal of Algebra}, 1(2):114--131, 1964.

\bibitem[FP23]{fernandez2023almost}
B.~Fern{\'a}ndez and S.~Penji{\'c}.
\newblock On (almost) 2-{$Y$}-homogeneous distance-biregular graphs.
\newblock {\em Bulletin of the Malaysian Mathematical Sciences Society},
  46(2):56, 2023.

\bibitem[GS70]{goethals1970strongly}
J.-M. Goethals and J.~J. Seidel.
\newblock Strongly regular graphs derived from combinatorial designs.
\newblock {\em Canadian Journal of Mathematics}, 22(3):597--614, 1970.

\bibitem[GST87]{distanceRegularised}
C.~Godsil and J.~Shawe-Taylor.
\newblock Distance-regularised graphs are distance-regular or
  distance-biregular.
\newblock {\em Journal of Combinatorial Theory, Series B}, 43(1):14--24, 1987.

\bibitem[HHL07]{huang2007class}
T.~Huang, L.~Huang, and M.~Lin.
\newblock On a class of strongly regular designs and quasi-semisymmetric
  designs.
\newblock In {\em Recent Developments in Algebra and Related Areas}, volume~8,
  pages 129--153. 2007.

\bibitem[Hig88]{higman1988strongly}
D.~Higman.
\newblock Strongly regular designs and coherent configurations of type [323].
\newblock {\em European Journal of Combinatorics}, 9(4):411--422, 1988.

\bibitem[HW68]{hall1968simple}
M.~Hall, Jr. and D.~Wales.
\newblock The simple group of order 604,800.
\newblock {\em Journal of Algebra}, 9(4):417--450, 1968.

\bibitem[Lat23]{phD}
S.~Lato.
\newblock {\em Distance-biregular graphs and orthogonal polynomials}.
\newblock PhD thesis, University of Waterloo, 2023.

\bibitem[Lat25]{soloCharacterization}
S.~Lato.
\newblock Polynomial characterizations of distance-biregular graphs.
\newblock {\em Journal of Graph Theory}, 109(3):282--293, 2025.

\bibitem[MST85]{biregularCage}
B.~Mohar and J.~Shawe-Taylor.
\newblock Distance-biregular graphs with 2-valent vertices and distance-regular
  line graphs.
\newblock {\em Journal of Combinatorial Theory, Series B}, 38(3):193--203,
  1985.

\bibitem[Neu80]{neumaier1980t12}
A.~Neumaier.
\newblock {$t \frac{1}{2}$}-designs.
\newblock {\em Journal of Combinatorial Theory, Series A}, 28(3):226--248,
  1980.

\bibitem[Neu82]{neumaier1982regular}
A.~Neumaier.
\newblock Regular sets and quasi-symmetric 2-designs.
\newblock In {\em Combinatorial Theory: Proceedings of a Conference Held at
  Schloss Rauischholzhausen, May 6--9, 1982}, pages 258--275. Springer, 1982.

\bibitem[Sec89]{secker1989feasibility}
P.~R. Secker.
\newblock {\em Feasibility algorithms for distance-biregular graphs}.
\newblock PhD thesis, University of Birmingham, 1989.

\bibitem[SS91]{shrikhande1991quasi}
M.~S. Shrikhande and S.~S. Sane.
\newblock {\em Quasi-symmetric designs}, volume 164.
\newblock Cambridge University Press, 1991.

\bibitem[ST85]{shaweTaylor}
J.~Shawe-Taylor.
\newblock {\em Regularity and transitivity in graphs}.
\newblock PhD thesis, University of London, 1985.

\bibitem[Tha95]{thas1995generalized}
J.~A. Thas.
\newblock Generalized polygons.
\newblock In {\em Handbook of Incidence Geometry}, pages 383--431. Elsevier,
  1995.

\bibitem[Tha07]{thas2007partial}
J.~A. Thas.
\newblock Partial geometries.
\newblock In {\em Handbook of Combinatorial Designs}, pages 557--561. Chapman
  \& Hall/CRC, 2007.

\bibitem[Tit59]{tits1959trialite}
J.~Tits.
\newblock Sur la trialit{\'e} et certains groupes qui s'en d{\'e}duisent.
\newblock {\em Publications Math{\'e}matiques de l'IH{\'E}S}, 2:13--60, 1959.

\bibitem[VM98]{van2012generalized}
H.~Van~Maldeghem.
\newblock {\em Generalized polygons}.
\newblock Springer Science \& Business Media, 1998.

\bibitem[Wit38]{witt}
E.~Witt.
\newblock {\"U}ber steinersche systeme.
\newblock {\em Abhandlungen aus dem Mathematischen Seminar der Universit{\"a}t
  Hamburg}, 12:265--275, 1938.

\bibitem[Yan81]{yanushka1981order}
A.~Yanushka.
\newblock On order in generalized polygons.
\newblock {\em Geometriae Dedicata}, 10(1--4):451--458, 1981.

\bibitem[Yos99]{MR1703601}
S.~Yoshiara.
\newblock A family of {$d$}-dimensional dual hyperovals in {PG}$(2d+1, 2)$.
\newblock {\em European Journal of Combinatorics}, 20(6):589--603, 1999.

\bibitem[Yos06]{Yoshiara_Pingree}
S.~Yoshiara.
\newblock Dimensional dual arcs---a survey.
\newblock In {\em Finite Geometries, Groups, and Computation}, pages 247--266.
  2006.

\end{thebibliography}

\end{document}